\documentclass[11pt]{amsart}

\usepackage{epic, eepic}
\usepackage{units}
\usepackage{amsthm, amsmath, amssymb,amscd}
\usepackage[all]{xy}
\usepackage{mathrsfs}
\usepackage{graphicx}
\usepackage{graphs}
\newcommand{\nf}[2]{\text{\large $\nicefrac{#1}{#2}$}}
\newtheorem*{claim}{Claim}
\newtheorem*{theorem}{Theorem}
\newtheorem{lemma}{Lemma}
\newtheorem*{remark}{Remark}
\newtheorem{example}{Example}
\newtheorem*{definition}{Definition}
\newtheorem*{corollary}{Corollary}
\newtheorem*{proposition}{Proposition}

\title{Hodge Spaces of Real Toric Varieties}
\author{Valerie Hower}
\address{Department of Mathematics, University of Georgia, Athens, GA 30602}
\email{vhower@math.uga.edu}

\begin{document}

\begin{abstract}

 We define the $\mathbb{Z}_2$ Hodge spaces $H_{pq}(\Sigma)$ of a fan $\Sigma$.
 If $\Sigma$ is the normal fan of a reflexive polytope $\Delta$ then we use polyhedral duality to compute the $\mathbb{Z}_2$ Hodge Spaces of $\Sigma$.  In particular, if the cones of dimension at most $e$ in the face fan $\Sigma^*$ of $\Delta$
 are smooth then we compute $H_{pq}(\Sigma)$ for $p<e-1$.   If $\Sigma^*$ is a smooth fan then we completely determine the spaces $H_{pq}(\Sigma)$ and we show $X_{\Sigma}$ is maximal, meaning
 that the sum of the $\mathbb{Z}_2$ Betti numbers of $X_{\Sigma}(\mathbb{R})$ is equal to the sum of the $\mathbb{Z}_2$ Betti numbers of $X_{\Sigma}(\mathbb{C})$.
\end{abstract}

\maketitle
\section{Introduction}
In this paper, we define and study the $\mathbb{Z}_2$ Hodge spaces of a fan $\Sigma$.  The $\mathbb{Z}_2$ Hodge spaces $H_{pq}(\Sigma)$ are indexed by pairs of integers $p,q$ with $0\leq q\leq p \leq d$,
 where $d=\mathrm{dim}\Sigma$.   When $\Sigma$ is a smooth fan, we have $H_{pq}(\Sigma)=0$ for $p\neq q$.  However, for $p>q$ the spaces $H_{pq}(\Sigma)$ are not generally well understood.  The terminology $\mathbb{Z}_2$ Hodge spaces is inspired by work of Brion [Bri] who considered similar spaces associated to a fan $\Sigma$. 

In Section \ref{hodgesec}, we introduce cosheaves on a fan $\Sigma$.  We also discuss cosheaf homology and some basic properties of cosheaves
on fans.  In later sections, we will also use
sheaves on $\Sigma$ and sheaf cohomology.  The definition and properties of sheaves on fans are similar to those of cosheaves
and are developed in [Bri \textsection 1.1].  Sheaves on fans are also studied in [Bre] and [Bar].
The main difference between our work with sheaves and that in [Bri] is that our sheaves are sheaves of $\mathbb{Z}_2$ vector
spaces. Our main interest lies in the cosheaf $\mathscr{E}$, which we define in Section \ref{esec}.  Geometrically, for $\sigma \in \Sigma$ the stalk $E_{\sigma}$ of the cosheaf $\mathscr{E}$ is the
compact real torus in the real orbit $O_{\sigma}(\mathbb{R})$ of the real toric variety $X_{\Sigma}(\mathbb{R})$.  The $\mathbb{Z}_2$ Hodge spaces $H_{pq}(\Sigma)$ are defined to be the homology groups $H_p(\wedge^q\mathscr{E}),$ where $\wedge^q\mathscr{E}$ is the $q$th exterior power of the cosheaf $\mathscr{E}$ on $\Sigma$.   

The $\mathbb{Z}_2$ Hodge spaces of $\Sigma$ are related to the topology of both the real and complex points of the toric variety $X_{\Sigma}$ in that $ H_{pq}(\Sigma)=\overline{E}_{p,q}^1=E_{p,q}^2,$ where $(\overline{E}^r, \overline{d}^r)$  and $(E^r, d^r)$ are two spectral sequences with 
 $$\xymatrix{\overline{E}^1_{p,q} \ar@{=>}[r]& H_p(X_{\Sigma}(\mathbb{R}), \mathbb{Z}_2) }
  \quad \mbox{and}\quad  \xymatrix{E^2_{p,q} \ar@{=>}[r]& H_{p+q}(X_{\Sigma}(\mathbb{C}), \mathbb{Z}_2) }. $$  The spectral sequences $(\overline{E}^r, \overline{d}^r)$ and $(E^r, d^r)$ are described in Section \ref{specsec}. 

 Section \ref{prelimsec} gives definitions and results which are needed for proofs in Section \ref{reflexsection}.  Section \ref{thomeosec} develops the notion of $T$-homeomorphic real toric varieties.  We develop a needed algebraic result in Section \ref{kossec} which states that two ``Koszul'' sequences of cosheaves are exact.  Section \ref{diagsec} gives an interpretation of the diagonal $\mathbb{Z}_2$ Hodge spaces of $\Sigma$ as the $\mathbb{Z}_2$ torus invariant Chow groups of $X_{\Sigma}$, and in Section \ref{rightsec} we compute the right-most column $H_{dq}(\Sigma)$.

 Section \ref{reflexsection} is devoted to establishing the following theorem.
 \begin{theorem}
 Let $\Delta$ be a reflexive polytope and $\Sigma^*$ the face fan of $\Delta$.  If the cones in $\Sigma^*$ of dimension $e$ are $\mathbb{Z}_2$ regular then
 \begin{enumerate}
 \item $H_{pq}(\Sigma)=0 \quad  \mbox{for}  \quad q<p<e-1$
 \item $H_{qq}(\Sigma)=\mathbb{Z}_2 \quad \mbox{for} \quad q< e-1$
 \end{enumerate}
 where $\Sigma$ is the normal fan of $\Delta$.
 \end{theorem}  An outline of Section \ref{reflexsection} is as follows.  Section \ref{corrsec} describes a correspondence between sheaves on $\Sigma^*$ and cosheaves on $\Sigma$, and in Section \ref{fcsection} we define the cosheaves $\widehat{\mathscr{F}}$, $\widehat{\mathscr{G}}$, and $\mathscr{C}$ on $\Sigma$.   We show that $0\longrightarrow \mathscr{C} \longrightarrow \widehat{\mathscr{G}} \longrightarrow \mathscr{E}^{\circ} \longrightarrow 0$ is a short exact sequence of cosheaves on $\Sigma$, where $H_p(\wedge^q\mathscr{E}^{\circ})=H_p(\wedge^q\mathscr{E})$ for $0\leq p \leq d-2$.  In Section \ref{vansec}, we use this short exact sequence and information about the vanishing of the homology groups $H_p(\wedge^k\widehat{\mathscr{G}}\, )$ (from Section \ref{gsec}) to obtain part (1) of the above theorem.    We prove part (2) of the theorem in Section \ref{chowsec}.

In Section \ref{collsec} we assume $\Sigma^*$ consists of $\mathbb{Z}_2$ regular cones.   We show that the spectral sequence $(\overline{E}^r, \overline{d}^r)$ for $X_{\Sigma}$ collapses at $\overline{E}^1.$ This gives that the sum of the $\mathbb{Z}_2$ Betti numbers for $X_{\Sigma}(\mathbb{R})$ is equal to the sum of the $\mathbb{Z}_2$ Betti numbers for $X_{\Sigma}(\mathbb{C})$.  Section \ref{exsec} contains two examples which illustrate our results. 

Some notation:  $N\cong \mathbb{Z}^d$ and $M=\mathrm{Hom}(N,\mathbb{Z})$ are dual lattices with dual pairing denoted $<\cdot,\cdot>$.  We write $\Sigma \subset N$ when $\Sigma \subset N\otimes \mathbb{R}$ is the normal fan of a $d$ dimensional lattice polytope $\Delta\subset M \otimes \mathbb{R}$ (written $\Delta \subset M$).  For $\sigma \in \Sigma$ we will abuse notation by writing $\sigma \cap N$ for the sublattice  $(\sigma \cap N) + (-\sigma \cap N)$.

All homology and cohomology groups will be with $\mathbb{Z}_2$ coefficients, unless
otherwise stated.

I would like to thank Matthias Franz and Tom Braden for helpful discussions at the June 2006 Toric Topology conference in Osaka, my advisor Clint McCrory for 
our ongoing discussions, and Bernd Sturmfels for introducing me to the beauty of real toric varieties.
\section{The $\mathbb{Z}_2$ Hodge spaces and the spectral sequence $\overline{E}^r$}
 \subsection{Cosheaves on a fan}\label{hodgesec}
A \emph{cosheaf} $\mathscr{F}$ of $\mathbb{Z}_2$ vector spaces on a fan $\Sigma$ is a collection of vector spaces $(F_{\sigma})_{\sigma\in \Sigma} $ over $\mathbb{Z}_2$  together with face restriction maps $\rho_{\tau, \sigma}: F_{\sigma}\longrightarrow F_{\tau}$ for $\sigma < \tau$ satisfying the following two conditions.
\begin{itemize}
\item If $ \sigma < \tau < \beta$ then $\rho_{\beta ,\tau}\rho_{\tau ,\sigma}= \rho_{\beta ,\sigma}$
\item $\rho_{\sigma,\sigma}$ is the identity map
\end{itemize}
We define \emph{cosheaf homology} groups as follows.    The chain groups $C_p(\mathscr{F})$ are the $\mathbb{Z}_2$ vector spaces defined by
$$C_p(\mathscr{F}):=\bigoplus_{\sigma \in \Sigma(d-p)}F_{\sigma}$$ and the boundary map $\partial_p: C_p(\mathscr{F})\longrightarrow C_{p-1}(\mathscr{F})$
is the direct sum of the maps \begin{equation}\label{nosign}\sum_{\sigma<\tau} \rho_{\tau,\sigma}: F_{\sigma} \longrightarrow \bigoplus_{{\small \tau \in \Sigma(d-p+1), \: \sigma<\tau}} F_{\tau}. \end{equation}
Note that if we were working with cosheaves of $k$ vector spaces with $\mathrm{char}k\neq 2$ then we would need to introduce
signs in (\ref{nosign}) to guarantee $\partial^2=0$.
We define $H_p(\mathscr{F})$ to be the $p$th homology group of the complex $(C_*(\mathscr{F}),\partial_*)$. \\
\indent Suppose $\mathscr{F}^1$ and $\mathscr{F}^2$ are two cosheaves on $\Sigma$.  A \emph{morphism of cosheaves} $\Phi: \mathscr{F}^1\longrightarrow \mathscr{F}^2$ is a collection of vector space homomorphisms $(\phi_{\sigma})_{\sigma\in \Sigma}$ with $\phi_{\sigma}: F^1_{\sigma} \longrightarrow F^2_{\sigma}$ such that if $\sigma < \tau<\beta$ the following diagram commutes.
$$\xymatrix{ F^1_{\sigma} \ar[d]^{\phi_{\sigma}} \ar[r]^{\rho^1_{\tau,\sigma}}&F^1_{\tau} \ar[d]^{\phi_{\tau}} \ar[r]^{\rho^1_{\beta,\tau}}&F^1_{\beta} \ar[d]^{\phi_{\beta}}    \\ F^2_{\sigma} \ar[r]^{\rho^2_{\tau, \sigma}}& F^2_{\tau} \ar[r]^{\rho^2_{\beta, \tau}} &F^2_{\beta} } $$
Hence, an \emph{exact sequence of cosheaves} $0\longrightarrow \mathscr{F}^1\longrightarrow \mathscr{F}^2 \longrightarrow \cdots \longrightarrow \mathscr{F}^n \longrightarrow 0$ is a collection of exact sequences of $\mathbb{Z}_2$ vector spaces for $\sigma \in \Sigma$
$$\xymatrix{0\ar[r] & F^1_{\sigma} \ar[r]^{\phi^1_{\sigma}}& F^2_{\sigma} \ar[r]^{\phi^2_{\sigma}}&\: \cdots \: \ar[r]^{\phi^{n-1}_{\sigma} }&F^n_{\sigma} \ar[r] &0 } $$ such that the maps $\phi^i_{\sigma}$ are natural with respect to face restriction for $1\leq i \leq n-1$.
In subsequent sections, we will use the fact that a short exact sequence of cosheaves
$$0\longrightarrow \mathscr{F}^1 \longrightarrow \mathscr{F}^2 \longrightarrow \mathscr{F}^3 \longrightarrow 0$$ induces a long exact
sequence on homology groups
$$\cdots \longrightarrow H_{p+1}(\mathscr{F}^3)\longrightarrow H_p(\mathscr{F}^1)\longrightarrow H_p(\mathscr{F}^2)\longrightarrow H_p(\mathscr{F}^3) \longrightarrow H_{p-1}(\mathscr{F}^1)\longrightarrow \cdots .$$
\subsection{The cosheaf $\mathscr{E}$ and the $\mathbb{Z}_2$ Hodge spaces $H_{pq}(\Sigma)$}\label{esec}
 We define the cosheaf $\mathscr{N}$ on $\Sigma$ by $$N_{\sigma}:=\nf{ \sigma \cap N}{ \sigma \cap 2N} \quad \mathrm{for} \quad
\sigma \in \Sigma$$ with the restriction map
$\rho_{\tau,\sigma}$ given by inclusion $$\xymatrix{\rho_{\tau,\sigma}: N_{\sigma} \ar[r]^{\subset \quad}& N_{\tau} \quad \mathrm{for} \quad \sigma < \tau }.$$
If $\mathrm{dim}\, \sigma=q$ then $N_{\sigma}$ is a rank $q$ vector space over $\mathbb{Z}_2$.  The cosheaf $\mathscr{E}$ is defined so that
$$E_{\sigma} :=N(\sigma) =\frac{\nf{N}{2N}}{N_{\sigma}}.$$
If $\sigma < \tau$ then the restriction map $\varpi_{\tau, \sigma}:N({\sigma}) \longrightarrow N({\tau})$
is induced from the identity on $\nf{N}{2N}$ which takes $N_{\sigma}$ to $N_{\tau}$.
Thus, the cosheaf $\mathscr{E}$ is the cokernel of the inclusion $\mathscr{N}\hookrightarrow \nf{N}{2N}$, where $\nf{N}{2N}$ is the constant
cosheaf assigning $\nf{N}{2N}$ to each cone in $\Sigma$. The $\mathbb{Z}_2$ \emph{Hodge spaces of }$\Sigma$ are defined to be the homology groups $H_{pq}(\Sigma):=H_p(\wedge^q\mathscr{E}),$ where $\wedge^q\mathscr{E}$ is the $q$th exterior power of
the cosheaf $\mathscr{E}$.
\begin{remark}
Since $\wedge^q E_{\sigma}$ is zero for $q>\mathrm{codim}\, \sigma$, we have $C_p(\wedge^q\mathscr{E})=0$ for $p<q$.
Hence the $\mathbb{Z}_2$ Hodge spaces $H_{pq}(\Sigma)$ are indexed by integers $p,q$ with $0\leq q\leq p \leq d$.
\end{remark}
\subsection{The spectral sequence $\overline{E}^r$}\label{specsec}
When  $X_{\Sigma}$ is projective, we have a natural cell structure on the real toric variety $X_{\Sigma}(\mathbb{R})$ given by the moment map.  We refer the reader to [Sot] for a discussion of real toric varieties and to [Gel \textsection 11.5 B] for a description of this cell structure.  The cells in $X_{\Sigma}(\mathbb{R})$ are of the form $(f,t)$
with $f< \Delta$ and  $t \in N(\sigma)$,  $\sigma \in \Sigma$ is dual to $f$.  
The boundary map $\partial$ for the chain complex $C_*(X_{\Sigma}(\mathbb{R}))$ is induced from the face restriction maps
$\varpi_{\tau, \sigma}: N(\sigma)\longrightarrow N(\tau)$ from Section \ref{hodgesec}.  It follows that
\begin{equation}\label{realchain}
C_j(X_{\Sigma}(\mathbb{R})) =\bigoplus_{\sigma \in \Sigma(d-j)}H_0(N(\sigma))
\end{equation}
where the boundary map $\partial_j$ is the direct sum of the maps
$$ \sum_{\sigma <\tau} (\varpi_{\tau, \sigma})_*:H_0(N(\sigma))\longrightarrow \bigoplus_{\tau \in \Sigma(d-j+1), \, \sigma <\tau} H_0(N(\tau)). $$

In [Bih] Bihan et al. want to understand the relationship between the topology of $X_{\Sigma}(\mathbb{R})$ and that of
$X_{\Sigma}(\mathbb{C})$.  The authors show that the chain complex $C_*(X_{\Sigma}(\mathbb{R}))$ is a filtered differential graded
vector space.  The associated spectral sequence $(\overline{E}^r, \overline{d}^r)$ converges to $H_*(X(\mathbb{R}))$ and is known to
collapse at $\overline{E}^1$ when $X_{\Sigma}$  is complete and has isolated singularities or when the dimension of $X_{\Sigma}$ is at most 3 [Bih].  Our notation is slightly different than in [Bih] and hence we briefly review their construction.\\
\indent Using (\ref{realchain}) we may specify a filtration on $C_*(X_{\Sigma}(\mathbb{R}))$ by giving a filtration of $H_0(N(\sigma))$, the $\mathbb{Z}_2$ group algebra of $N(\sigma)$.  We use the \emph{augmentation} homomorphism $\epsilon_{\sigma}$.
\begin{eqnarray*}
\epsilon_{\sigma}: H_0(N(\sigma))& \longrightarrow &\mathbb{Z}_2 \\  \sum_{n_i\in \mathbb{Z}_2, \, g_i \in  N(\sigma)}n_ig_i &\longmapsto & \sum n_i
\end{eqnarray*}
We define $I_{\sigma}$, an ideal in $H_0(N(\sigma))$, via $I_{\sigma}:=\mathrm{ker}\, \epsilon_{\sigma}$.   This gives a filtration of $H_0(N(\sigma))$
$$
0=I_{\sigma}^{j+1} \subset I_{\sigma}^j \subset \cdots \subset I_{\sigma}^2 \subset I_{\sigma} \subset I_{\sigma}^0= H_0(N(\sigma))$$
where $j=\mathrm{rank}\,N(\sigma)=\mathrm{codim}\,\sigma$.
We reindex by setting $J_{\sigma}^p=I_{\sigma}^{d-p}$ so that
\begin{equation}\label{filter}
0=J_{\sigma}^{d-j-1} \subset J_{\sigma}^{d-j} \subset \cdots \subset J_{\sigma}^{d-2} \subset J_{\sigma}^{d-1} \subset J_{\sigma}^{d}= H_0(N(\sigma))
\end{equation}
is an increasing filtration of $H_0(N(\sigma))$.
\begin{lemma} \label{filtlem}
The filtrations of $H_0(N(\sigma))$ of the form (\ref{filter}) for $\sigma \in \Sigma$ determine an increasing filtration of $F$ of $C_*(X_{\Sigma}(\mathbb{R}))$ and we have the following. $$\partial_j:F_q C_j(X_{\Sigma}(\mathbb{R}))\longrightarrow F_q C_{j-1}(X_{\Sigma}(\mathbb{R}))$$
\end{lemma}
To prove Lemma \ref{filtlem} note that for $\sigma < \tau$ in $\Sigma$ the map $(\varpi_{\tau,\sigma})_*$ commutes with the augmentation homomorphisms.
$$\xymatrix{ H_0(N(\sigma)) \ar[rr]^{(\varpi_{\tau,\sigma})_*} \ar[dr]_{\epsilon_{\sigma}}& & H_0(N(\tau)) \ar[dl]^{\epsilon_{\tau}} \\ & \mathbb{Z}_2&}$$
Thus we have $(\varpi_{\tau,\sigma})_* ( I_{\sigma} ) \subset I_{\tau} $ and $(\varpi_{\tau,\sigma})_* (I_{\sigma}^k) \subset I_{\tau}^k$.  Lemma \ref{filtlem} follows.  We will denote $\underline{\varpi}^k_{\tau,\sigma}$ for the induced map on the quotient
$$\underline{\varpi}^k_{\tau,\sigma}: \frac{I_{\sigma}^k}{I_{\sigma}^{k+1}} \longrightarrow \frac{I_{\tau}^k}{I_{\tau}^{k+1}}.$$
\indent We now use the spectral sequence of a
filtered module as developed in [Mac \textsection XI.3] to define the spectral sequence $(\widetilde{E}^r, \widetilde{d}^r)$.
From Theorem 3.1 in [Mac] we have $$\widetilde{E}^1_{p,q}=H_{p+q}\left(\frac{F_p\,C_*(X_{\Sigma}(\mathbb{R}))}{F_{p-1}\,C_*(X_{\Sigma}(\mathbb{R}))}\right)$$
and $$\widetilde{E}^1_{p,q} \Longrightarrow H_{p+q}(X_{\Sigma}(\mathbb{R})).$$
\begin{lemma}\label{overlinelemma}
For each $p,q$ we have
$$ \widetilde{E}^1_{p,q} \cong H_{p+q, \, d-p}(\Sigma)$$ where $H_{p+q, \, d-p}(\Sigma)$ is the $\mathbb{Z}_2$ Hodge space of $\Sigma$ defined in Section \ref{hodgesec}.
\end{lemma}
To prove Lemma \ref{overlinelemma} consider the following.
\begin{eqnarray*}
& \text{{ \large $\frac{F_p\,C_{p+q}(X_{\Sigma}(\mathbb{R}))}{F_{p-1}\,C_{p+q}(X_{\Sigma}(\mathbb{R}))}$ }}
& =\bigoplus_{\sigma \in \Sigma(d-(p+q))}\frac{J_{\sigma}^p}{J_{\sigma}^{p-1}} \\
&& =\bigoplus_{\sigma \in \Sigma(d-(p+q))}\frac{I_{\sigma}^{d-p}}{I_{\sigma}^{d-(p-1)}}\\
\end{eqnarray*}
The following claim is Proposition 6.1 in [Bih].
\begin{claim} For each $\sigma \in \Sigma$ we have
$$\frac{I_{\sigma}^q}{I_{\sigma}^{q+1}} \cong \wedge^q N(\sigma)$$ canonically as $\mathbb{Z}_2$ vector spaces.
\end{claim}
 This yields
\begin{eqnarray*}
& \text{{\large $\frac{F_p\,C_{p+q}(X_{\Sigma}(\mathbb{R}))}{F_{p-1}\,C_{p+q}(X_{\Sigma}(\mathbb{R}))}$}}
& =\bigoplus_{\sigma \in \Sigma(d-(p+q))} \wedge^{d-p}N(\sigma) \\
&& =C_{p+q}(\wedge^{d-p} \mathcal{E}). \end{eqnarray*}
To see the lemma we need only compare boundary maps.  The map
$$\widetilde{d}^0_{p,q}:  \frac{F_p\,C_{p+q}(X_{\Sigma}(\mathbb{R}))}{F_{p-1}\,C_{p+q}(X_{\Sigma}(\mathbb{R}))} \longrightarrow
 \frac{F_{p}\,C_{p+q-1}(X_{\Sigma}(\mathbb{R}))}{F_{p-1}\,C_{p+q-1}(X_{\Sigma}(\mathbb{R}))} $$ is given by the
collection
$$ \underline{\varpi}^{d-p}_{\tau,\sigma} : \, \frac{I_{\sigma}^{d-p}}{I_{\sigma}^{d-p+1}} \longrightarrow \frac{I_{\tau}^{d-p}}{I_{\tau}^{d-p+1}} $$
of maps for $\sigma < \tau$ in $\Sigma$ with $\mathrm{dim}\, \sigma=d-(p+q)$.
By construction, this is induced from the map $\varpi_{\tau,\sigma}$.  The cosheaf differential for $\wedge^{d-p}\mathscr{E}$ is given by the collection of face restriction maps
$$ \wedge^{d-p}N(\sigma) \longrightarrow \wedge^{d-p}N(\tau) \quad \mathrm{for} \quad \sigma < \tau$$
which are also induced from the maps $\varpi_{\tau,\sigma}$, and hence the two maps are equal.  Taking homology we arrive at Lemma \ref{overlinelemma}.\\
\indent As $(\widetilde{E}^r, \widetilde{d}^r)$ is the spectral sequence associated to the filtered complex of cellular chains in $X_{\Sigma}(\mathbb{R})$, we have
$$\widetilde{E}^{\infty}_{p,q}= \frac{F_p\,H_{p+q}(X_{\Sigma}(\mathbb{R}))}{F_{p-1}\,H_{p+q}(X_{\Sigma}(\mathbb{R}))}$$ with
{\small $$0=F_{d-(p+q)-1}H_{p+q}(X_{\Sigma}(\mathbb{R})) \subset F_{d-(p+q)}H_{p+q}(X_{\Sigma}(\mathbb{R})) \subset  \cdots \subset F_{d}H_{p+q}(X_{\Sigma}(\mathbb{R})) \subset H_{p+q}(X_{\Sigma}(\mathbb{R}))$$}
the bounded filtration of $H_{p+q}(X_{\Sigma}(\mathbb{R}))$.  Since $\widetilde{d}^r_{d-(p+q),2(p+q)-d}=0$ for all $r$, we have a sequence of surjections
$$\xymatrix{\widetilde{E}^1_{d-(p+q),2(p+q)-d} \ar@{->>}[r]& \widetilde{E}^2_{d-(p+q),2(p+q)-d} \ar@{->>}[r]&  \widetilde{E}^3_{d-(p+q),2(p+q)-d} \ar@{->>}[r]&  \cdots \ar@{->>}[r]& \widetilde{E}^{\infty}_{d-(p+q),2(p+q)-d} }$$
which are called edge homomorphisms in \textsection 11.1 of [Mac].
Moreover, as $$\widetilde{E}^{\infty}_{d-(p+q),2(p+q)-d} =\frac{F_{d-(p+q)}H_{p+q}(X_{\Sigma}(\mathbb{R}))}{0},$$ we can compose the edge homomorphisms with inclusion to obtain a natural homomorphism
\begin{equation}\label{edge}
\widetilde{E}^1_{d-(p+q),2(p+q)-d} \longrightarrow H_{p+q}(X_{\Sigma}(\mathbb{R})).
\end{equation}
\indent At this point, we introduce reindexing which we will use in the sequel.  We rotate the $\widetilde{E}^1$ term clockwise by $90^{\circ}$ and then shear so that the
diagonal lines $p+q=c$ become vertical lines $p=c$.  We write $(\overline{E}^r, \overline{d}^r)$ for the spectral sequence obtained after
reindexing.  
\begin{figure}[htbp]\vspace{-.55in}
\begin{picture}(8,4) \put(-3.8,1.4){{\footnotesize $q$}}
\put(-3.5,0){ {\footnotesize $ 
 \begin{array}{|ccccc} H_{44}(\Sigma ) && & &
\\  & H_{43}(\Sigma)& && \\  &  H_{33}(\Sigma ) & H_{42}(\Sigma )&&\\  && H_{32}(\Sigma ) & H_{41}(\Sigma ) &  \\
&& H_{22}(\Sigma )& H_{31}(\Sigma ) &H_{40}(\Sigma)   \\\hline  \vspace{-.3cm} \\&p&& H_{21}(\Sigma )& H_{30}(\Sigma ) \\ & && H_{11}(\Sigma ) & H_{20}(\Sigma ) \\
&&&& H_{10}(\Sigma ) \\ &&&& H_{00}(\Sigma)   \\ \end{array} $ }}
\put(4.7,1){{\footnotesize $q$}}
\put(5,0){ {\footnotesize $ 
 \begin{array}{|ccccc}  &&&& H_{44}(\Sigma ) \\&&&  H_{33}(\Sigma ) & H_{43}(\Sigma )
\\ && H_{22}(\Sigma ) & H_{32}(\Sigma ) & H_{42}(\Sigma ) \\ & H_{11}(\Sigma ) & H_{21}(\Sigma ) & H_{31}(\Sigma )
& H_{41}(\Sigma )    \\ H_{00}(\Sigma )& H_{10}(\Sigma ) & H_{20}(\Sigma ) & H_{30}(\Sigma )
& H_{40}(\Sigma )  \\ \hline \vspace{-0.3cm} \\ &p&  && \\ \end{array} $
}}
\end{picture}
\vspace{2.3cm}
 \caption{The terms $\widetilde{E}^1_{p,q}$ (left) and  $\overline{E}^1_{p,q}$ (right)}\label{tildebar}
\end{figure} 
 The total grading of $(\overline{E}^r, \overline{d}^r)$ is $p$, and Lemma \ref{overlinelemma} gives the identity
$\overline{E}^1_{p,q}=H_{pq}(\Sigma).$  Figure \ref{tildebar} shows the terms $\widetilde{E}^1_{p,q}$  and  $\overline{E}^1_{p,q}$ for a 4 dimensional fan $\Sigma$.
 
For the spectral sequence $(\overline{E}^r,\overline{d}^r)$ the boundary map $\overline{d}^r$ satisfies
$\overline{d}^r_{p,q}:\overline{E}^r_{p,q} \longrightarrow \overline{E}^r_{p-1, q+r}.$ We show the
 first five boundary maps in Figure \ref{diffpic}. \vspace{.2cm}
  \begin{figure}[htbp]    \begin{center} \setlength{\unitlength}{1.1cm}     \begin{graph}(1.5,3.5)(0,.9)
  \graphnodesize{.1} \graphlinewidth{.02}\grapharrowlength{.16}\grapharrowwidth{.5}
\roundnode{1}(0.2,0.2)  \roundnode{2a}(1.2,1.2)
\roundnode{2}(0.2,1.2) \roundnode{3a}(1.2,2.2)
\roundnode{3}(0.2,2.2) \roundnode{4a}(1.2,3.2)
\roundnode{4}(0.2,3.2) \roundnode{5a}(1.2,4.2)
\roundnode{5}(0.2,4.2)
\roundnode{6}(1.2,0.2)\diredge{6}{5}
\diredge{6}{1}
\diredge{6}{2}
\diredge{6}{4}
\diredge{6}{3}
\freetext(.1,.45){{\footnotesize $\overline{d}^0$}} \freetext(.1,4.45){{\footnotesize $\overline{d}^4$}}
\freetext(.1,1.45){{\footnotesize $\overline{d}^1$}} \nodetext{6}(0,-.2){ {\scriptsize $(p,q)$}}
\freetext(.1,3.45){{\footnotesize $\overline{d}^3$}} \nodetext{1}(-.1,-.2){ {\scriptsize $(p-1,q)$}}
\freetext(.1,2.45){{\footnotesize $\overline{d}^2$}}
   \end{graph} \vspace{1cm} \caption{The differentials $\, \overline{d}^0, \overline{d}^1, \overline{d}^2, \overline{d}^3,$ and $\overline{d}^4$ }
\label{diffpic} \end{center}  \end{figure}
Next, we point out that after reindexing the natural homomorphism (\ref{edge}) becomes a map which we denote $f_*^{\mathbb{R}}$ and will use in future sections where
\begin{equation}\label{fmap}
f_{q}^{\mathbb{R}}: H_{qq}(\Sigma)=\overline{E}^1_{q,q}\longrightarrow H_q(X_{\Sigma}(\mathbb{R})).
\end{equation}
\begin{lemma} \label{ssequal}
The spectral sequence $(\overline{E}^r, \overline{d}^r)$ satisfies
$$\overline{E}^1_{p,q}\cong E^2_{p,q} \quad \mathrm{for} \quad 0\leq p,q \leq d$$ where $(E^r,d^r)$ is the $\mathbb{Z}_2$
Leray spectral sequence for the map $\mu_{\mathbb{C}}$.
\end{lemma}
We refer the reader to [Mac \textsection XI.7] or [McC \textsection 5] for construction of the Leray spectral sequence and to [Ful1 \textsection 4.2] for properties of the moment map $\mu_{\mathbb{C}}$.
We will show $E^2_{p,q} \cong H_{pq}(\Sigma)$ which will prove the lemma.  Suppose $\sigma \in \Sigma$ is dual to $f<\Delta$.  From  [Ful1 \textsection 4.2], if $z\in \mathrm{int}\,f$ then
$$\mu_{\mathbb{C}}^{-1}\{z\}= \frac{\nf{(N\otimes \mathbb{R})}{\sigma}}
{\nf{N}{\sigma\cap N}}\cong (S^1)^{p} \quad \quad p=\mathrm{dim}\,f$$
Next, we note that $H_q((S^1)^{p}, \mathbb{Z}_2) = \wedge^q\left(H_1((S^1)^{p}, \mathbb{Z}_2)\right)$ and
\begin{eqnarray*}
&H_1\left(  \frac{\nf{(N\otimes \mathbb{R})}{\sigma}}
{\nf{N}{\sigma\cap N}} , \mathbb{Z}_2\right)& =H_1\left(  \frac{\nf{(N\otimes \mathbb{R})}{\sigma}}
{\nf{N}{\sigma\cap N}} , \mathbb{Z}\right)\otimes \mathbb{Z}_2 \\
&&=(\nf{N}{\sigma\cap N})\otimes \mathbb{Z}_2 \\
&&=\frac{\nf{N}{2N}}{\nf{\sigma\cap N}{\sigma\cap 2N}} \\
&& =N(\sigma)
\end{eqnarray*}
where each equality above is canonical.
This gives the following.\begin{eqnarray*} &E^1_{p,q} &=\bigoplus_{\sigma \in \Sigma(d-p)}H_q\left( \frac{\nf{(N\otimes \mathbb{R})}{\sigma}}
{\nf{N}{\sigma\cap N}}, \mathbb{Z}_2\right) \\ && =\bigoplus_{\sigma \in \Sigma(d-p)}\wedge^q N(\sigma)
\end{eqnarray*}  Using Proposition 5.1 of [Bih], the differential $d^1_{*,q}: E^1_{*,q} \longrightarrow E^1_{*-1,q}$
is equivalent to the boundary map for the cosheaf $\wedge^q \mathscr{E}$ and hence we have Lemma \ref{ssequal}. \\
\indent We now have a way to compare the topology of the real and complex points of a toric variety.
The Smith-Thom inequality states the sum of the $\mathbb{Z}_2$ Betti numbers for $X_{\Sigma}(\mathbb{R})$ is less than or equal to the sum of the $\mathbb{Z}_2$ Betti numbers for $X_{\Sigma}(\mathbb{C})$.  We say $X_{\Sigma}$ is \emph{maximal} if equality is obtained.
Using the the Smith-Thom inequality and Lemma \ref{ssequal} we have the following diagram.
\begin{equation}\label{smith}
\xymatrix{\sum \mathrm{rank}(\overline{E}^1_{p,q})\ar@{}[r]|{=}\ar@{}[d]|{\vee |}& \sum \mathrm{rank}(E^2_{p,q})\ar@{}[d]|{\vee|}\\
 \sum b_i(X(\mathbb{R}))\ar@{}[r]|{\leq} & \sum b_j(X(\mathbb{C}))
} \end{equation}
If the spectral sequence $\overline{E}^r$ collapses at $\overline{E}^1$ then the left vertical inequality in (\ref{smith}) is equality.
This forces both the right vertical inequality and the lower horizontal inequality to both be equalities.  In this case we obtain $X_{\Sigma}$ is maximal.

\section{Preliminary results}\label{prelimsec}
\subsection{$T$-homeomorphic versus algebraically isomorphic} \label{thomeosec}
Suppose $\Sigma_1\subset N_1$ and $\Sigma_2 \subset N_2$ are two fans with $N_1\cong \mathbb{Z}^d$ and $N_2\cong \mathbb{Z}^n$.  Following [Oda \textsection 1.5],
let $A$ be a $d \times n$ integer valued matrix such that for each cone $\sigma_1 \in \Sigma_1$ there is a cone
$\sigma_2\in \Sigma_2$ with $A\sigma_1 \subset \sigma_2$.  Then, $A$ gives a \emph{map of fans} and determines a map of toric varieties.
$$\widetilde{A}:X_{\Sigma_1}\longrightarrow X_{\Sigma_2}$$ The restriction of $\widetilde{A}$ to the dense torus $T_1\subset X_{\Sigma_1}$
coincides with the homomorphism of tori
$$\xymatrix{T_1=N_1\otimes \mathbb{C}^* \ar[r]^{A\otimes 1} &N_2\otimes \mathbb{C}^*=T_2.
}$$
When $d=n$, the map $\widetilde{A}$ is an \emph{algebraic isomorphism} of $X_{\Sigma_1}$ and $X_{\Sigma_2}$ provided the cones in $\Sigma_2$ are all of the form $A\sigma$ for $\sigma\in \Sigma_1$ and $|\mathrm{det}A|=1$. \\
 \indent Two real toric varieties may be homeomorphic even though they are not isomorphic as algebraic varieties.  This motivates the notion of $T$-homeomorphic real toric varieties.
Let's assume $A$ is a $d\times d$ integer valued matrix with $\mathrm{det}A$ an odd integer.  We consider $A$ as a map between the $d$ dimensional lattices $N_1$ and $N_2$.  Then we have $A(2N_1) \subset 2N_2$ and hence $A$ determines a map
$$A^{\prime}: \nf{N_1}{2N_1} \longrightarrow \nf{N_2}{2N_2}.$$
Since the determinant of $A$ is odd, the preimage $A^{-1}\{2N_2\}$ is contained in $2N_1$.  Thus $A^{\prime}$ is an injective map between $\mathbb{Z}_2$ vector spaces of the same rank yielding $A^{\prime}$ is an isomorphism of vector spaces.
\begin{definition}
Suppose the $d\times d$ matrix $A$ determines a map of fans and $\mathrm{det}A$ is an odd integer.  Then $A$ induces a \emph{$T$-homeomorphism} $\xymatrix{ \widetilde{A}:X_{\Sigma_1}(\mathbb{R}) \ar[r]^{\cong} & X_{\Sigma_2}(\mathbb{R}) }$ provided  \begin{enumerate}
\item For each $\sigma_1\in \Sigma_1$ there exists $\sigma_2 \in \Sigma_2$ such that the restriction $$\bigl.A^{\prime}\bigr|_{N_{\sigma_1}}:N_{\sigma_1} \longrightarrow N_{\sigma_2}$$ is a vector space isomorphism.
\item For each $\sigma_2 \in \Sigma_2$, $N_{\sigma_2}=A^{\prime}(N_{\sigma_1})$ for some $\sigma_1\in \Sigma_1$.
\end{enumerate}  We say $X_{\Sigma_1}(\mathbb{R})$ and $X_{\Sigma_2}(\mathbb{R})$ are \emph{$T$-homeomorphic}.
\end{definition}
\subsection{Exactness of ``Koszul'' sequences}\label{kossec}
For a $\mathbb{Z}_2$ vector space $E$ of rank $r$, we have the Koszul complex as in [McC p. 259] or [Lan p. 861]
\begin{equation}\label{mccleary}
0\longrightarrow \wedge^rE\otimes SE\longrightarrow \wedge^{r-1}E\otimes SE \longrightarrow \cdots \longrightarrow \wedge^1E \otimes SE  \longrightarrow \wedge^0E\otimes SE\longrightarrow 0
\end{equation}
with boundary map $\, \partial_p:\wedge^pE \otimes SE \longrightarrow \wedge^{p-1}E \otimes SE$ given by
\begin{equation}\label{kosdiff}
\partial_p(x_i\wedge \cdots \wedge x_p\otimes y)= \sum_{i=1}^p x_1\wedge \cdots \wedge \widehat{x_i}\wedge \cdots \wedge x_p\otimes x_iy
\end{equation}
The complex (\ref{mccleary}) has $H_0=\mathbb{Z}_2$ and $H_i=0$ for  $i>0$, and hence
\begin{equation}\label{mcexact}
0\longrightarrow \wedge^rE\otimes SE\longrightarrow
\wedge^{r-1}E\otimes SE \longrightarrow \cdots \longrightarrow \wedge^1E \otimes SE  \longrightarrow \wedge^0E\otimes SE\longrightarrow \mathbb{Z}_2\longrightarrow 0
\end{equation}
is exact.  Decomposing (\ref{mcexact}) into its graded pieces, we obtain for $j\ge 1$ the following exact sequence of $\mathbb{Z}_2$ vector spaces.
\begin{equation}\label{strandj}
0\longrightarrow \wedge^jE\otimes S^0E\longrightarrow \wedge^{j-1}E\otimes S^1E \longrightarrow \cdots \longrightarrow
\wedge^1E \otimes S^{j-1}E  \longrightarrow \wedge^0E\otimes S^jE\longrightarrow 0
\end{equation}
Next, let $G$ be a $\mathbb{Z}_2$ vector space of rank $p$ and we consider the following for $\, k\ge 1$
{\footnotesize \begin{equation}\label{koz1}
0\longrightarrow (\bigoplus_{i+j=k}\wedge^iG\otimes \wedge^jE)\otimes S^0E\longrightarrow (\bigoplus_{l+t=k-1}\wedge^lG \otimes \wedge^{t}E)
\otimes S^1E \longrightarrow \cdots \longrightarrow \wedge^0G
\otimes \wedge^0E\otimes S^kE\longrightarrow0
\end{equation} }
where the boundary map is given by
\begin{eqnarray*}
&\wedge^iG\otimes \wedge^jE\otimes S^lE \longrightarrow \wedge^iG\otimes \wedge^{j-1}E \otimes S^{l+1}E \\
&\alpha \otimes \beta \otimes \gamma \longmapsto \alpha \otimes \partial_j(\beta \otimes \gamma)
\end{eqnarray*}
We use the sequence (\ref{koz1}) to construct the exact sequences (\ref{ekoz1}).
We use the fact that $$(\bigoplus_{i+j=q}\wedge^iG\otimes \wedge^{j}E)\otimes S^{k-q}E \cong
\bigoplus_{i+j=q}(\wedge^iG\otimes \wedge^jE\otimes S^{k-q}E )$$ and that the boundary map in (\ref{koz1}) is a direct sum.
On one summand, we obtain
\begin{eqnarray*}
&& \mathrm{ker}(\wedge^iG\otimes \wedge^jE\otimes S^{k-q}E \longrightarrow \wedge^iG\otimes \wedge^{j-1}E\otimes S^{k-q+1}E)\\
&& =\wedge^iG \otimes \mathrm{ker}( \wedge^jE\otimes S^{k-q}E \longrightarrow \wedge^{j-1}E\otimes S^{k-q+1}E) \\
&& \cong \left\{ \begin{array}{ll}\wedge^iG \otimes \mathrm{im}(\wedge^{j+1}E
\otimes S^{k-q-1}E \longrightarrow  \wedge^{j}E\otimes S^{k-q}E) & \mathrm{if} \; \; k-q\neq 0\\
0 & \mathrm{if} \; \; k-q=0 \; \;  \mathrm{and} \; \; j\neq 0 \\
\wedge^kG\otimes \wedge^0E\otimes S^0E & \mathrm{if} \; \; k-q=0 \; \; \mathrm{and} \; \; j=0. \\
\end{array} \right. \end{eqnarray*}
Since
\begin{eqnarray*} &&
\wedge^iG \otimes \mathrm{im}(\wedge^{j+1}E
\otimes S^{k-q-1}E \longrightarrow  \wedge^{j}E\otimes S^{k-q}E)\\
 && \cong \mathrm{im}(\wedge^iG\otimes\wedge^{j+1}E
\otimes S^{k-q-1}E \longrightarrow \wedge^iG\otimes  \wedge^{j}E\otimes S^{k-q}E) ,
\end{eqnarray*}
we have exactness at each piece except for the far left and
\begin{eqnarray*}
&&\mathrm{ker}(\bigoplus_{i+j=k}\wedge^iG\otimes \wedge^jE\otimes S^0E \longrightarrow \bigoplus_{i+j-1=k-1}\wedge^iG\otimes \wedge^{j-1}E \otimes S^1E) \\
&&=\wedge^kG\otimes \wedge^0E\otimes S^0E \\
&&\cong \wedge^kG .
\end{eqnarray*}
Hence, {\small
\begin{equation}\label{ekoz1}
0\longrightarrow \wedge^kG\longrightarrow (\bigoplus_{i+j=k}\wedge^iG\otimes \wedge^jE)\otimes S^0E\longrightarrow \cdots \longrightarrow \wedge^0G
\otimes \wedge^0E\otimes S^kE\longrightarrow0
\end{equation} }
is an exact sequence of $\mathbb{Z}_2$ vector spaces.

Next, suppose we have a short exact sequence of $\mathbb{Z}_2$ vector spaces
\begin{equation}\label{ses}
 \xymatrix{ 0\ar[r] &G \ar[r]^{\Phi}&F\ar[r]^{\Psi}&E\ar[r] &0 .\\
} \end{equation}
 We can form
\begin{equation}\label{kosa}
0\longrightarrow \wedge^kG \longrightarrow \wedge^kF \longrightarrow \wedge^{k-1}F\otimes S^1E \longrightarrow \cdots \longrightarrow \wedge^1F\otimes S^{k-1}E\longrightarrow S^kE\longrightarrow 0
\end{equation}
with the boundary map
\begin{eqnarray*} &&
\wedge^lF\otimes S^{k-l}E \longrightarrow \wedge^{l-1}F \otimes S^{k-l+1}E\\
&& x_1\wedge \cdots \wedge x_l\otimes a \longmapsto \sum_{i=1}^l x_i\wedge \cdots \wedge \widehat{x_i}\wedge \cdots \wedge x_l \otimes \Psi(x_i)\,a
\end{eqnarray*}
where  $$\sum_{i=1}^l x_i\wedge \cdots \wedge \widehat{x_i}\wedge \cdots \wedge x_l \otimes \Psi(x_i)\,a= \sum_{x_i\in \mathrm{coker}\Phi} x_i\wedge \cdots \wedge \widehat{x_i}\wedge \cdots \wedge x_l \otimes \Psi(x_i)\,a$$
because $x_i\in \Phi(G)  \iff \Psi(x_i)=0$.
Note that the boundary map in (\ref{kosa}) is equivalent to the boundary map in (\ref{ekoz1}), where we have used that $\displaystyle \bigoplus_{i+j=l}\wedge^iE\otimes \wedge^jG \cong \wedge^lF$.  This shows the sequence (\ref{kosa}) is exact.  Moreover, the boundary map in (\ref{kosa}) is independent of the choice of splitting of the short exact sequence (\ref{ses}). \\

\indent Next, we will obtain another exact sequence.  If $0\longrightarrow E\longrightarrow H\longrightarrow G\longrightarrow 0$ is exact then we have  $0\longrightarrow G^*\longrightarrow H^*\longrightarrow E^*\longrightarrow 0$ is an exact sequence of $\mathbb{Z}_2$ vector spaces.  In this case the exact sequence (\ref{kosa}) gives
{\small \begin{equation}\label{kosad}
0\longrightarrow \wedge^kG^* \longrightarrow \wedge^kH^* \longrightarrow \wedge^{k-1}H^*\otimes S^1E^* \longrightarrow \cdots \longrightarrow \wedge^1H^*\otimes S^{k-1}E^*\longrightarrow S^kE^*\longrightarrow 0.
\end{equation} }
Applying $\mathrm{Hom}(\relbar,\mathbb{Z}_2)$ to (\ref{kosad}) we obtain the exact sequence
\begin{equation}\label{kosc}
0\longrightarrow S^kE\longrightarrow \wedge^1H\otimes S^{k-1}E\longrightarrow \cdots \longrightarrow \wedge^{k-1}H\otimes S^1E\longrightarrow \wedge^kH\longrightarrow \wedge^kG \longrightarrow0
\end{equation}
where we have used that
\begin{eqnarray*}
(\wedge^jH^* \otimes S^{k-j}E^*)^* &&\cong (\wedge^jH^*)^*\otimes(S^{k-j}E^*)^* \\
&& \cong \wedge^jH\otimes S^{k-j}E
\end{eqnarray*} as vector spaces over $\mathbb{Z}_2$.  Moreover, the boundary maps in (\ref{kosc}) do not depend on a choice of splitting of the short exact sequence  $0\longrightarrow E\longrightarrow H\longrightarrow G\longrightarrow 0$ because the maps in (\ref{kosad}) do not depend on a splitting of $0\longrightarrow G^*\longrightarrow H^*\longrightarrow E^*\longrightarrow 0$. \\

We conclude this section with our interest in the exact sequences (\ref{kosa}) and (\ref{kosc}). 
\begin{proposition}
Suppose
\begin{equation}\label{sheafexact}
0\longrightarrow \mathscr{A}\longrightarrow \mathscr{B} \longrightarrow \mathscr{C} \longrightarrow 0
\end{equation} is an exact sequence of sheaves (or cosheaves) on a fan $\Sigma$.
Then the following sequences are exact for $k\ge 1$.
{ \small \begin{eqnarray}
&0\longrightarrow\wedge^k\mathscr{A} \longrightarrow\wedge^k\mathscr{B} \longrightarrow \wedge^{k-1}\mathscr{B}\otimes S^1\mathscr{C}
\longrightarrow \cdots \longrightarrow\wedge^1\mathscr{B}\otimes S^{k-1}\mathscr{C}\longrightarrow S^k\mathscr{C}\longrightarrow 0 \label{kozash} \\
&0\longrightarrow S^k\mathscr{A}\longrightarrow\wedge^1\mathscr{B}\otimes S^{k-1}\mathscr{A}\longrightarrow\cdots \longrightarrow
\wedge^{k-1}\mathscr{B}\otimes S^1\mathscr{A}\longrightarrow\wedge^k\mathscr{B}\longrightarrow\wedge^k\mathscr{C} \longrightarrow 0 \label{koscsh}
\end{eqnarray} }
\end{proposition}
For each $\sigma\in\Sigma$, we may apply (\ref{kosa}) and (\ref{kosc}) to the short exact sequence of
$\mathbb{Z}_2$ vector spaces given by the stalks $0\longrightarrow A_{\sigma}\longrightarrow B_{\sigma} \longrightarrow C_{\sigma} \longrightarrow 0.$
Since each boundary map in (\ref{kosa}) and (\ref{kosc}) is independent of a choice of splitting
of the short exact sequence $0\longrightarrow A_{\sigma}\longrightarrow B_{\sigma} \longrightarrow C_{\sigma} \longrightarrow 0$ and is induced from one of the the sheaf (or cosheaf) maps in  (\ref{sheafexact}),
the collection of maps on stalks gives rise to a sheaf (or cosheaf) homomorphism.
Thus, we obtain the exact sequences of sheaves (or cosheaves) in the proposition.


\subsection{The diagonal entries $H_{qq}(\Sigma)$}\label{diagsec}
In this section we work with the $\mathbb{Z}_2$ torus invariant Chow groups of $X_{\Sigma}$.  The integral
torus invariant Chow groups are discussed in  [Ful1 \textsection 5.1] or [Ful2].  When considering the $\mathbb{Z}_2$ torus invariant Chow groups, we define
$A^T_q(X_{\Sigma}):=\nf{Z^T_q(X_{\Sigma})}{R^T_{q+1}(X_{\Sigma})}$, where $Z^T_q(X_{\Sigma})$ is the $\mathbb{Z}_2$ vector space
generated by cycles $[V(\tau)]$ for $\tau \in \Sigma$ of codimension $q$.  The relations are generated by torus invariant divisors of the following form.  For $\sigma \in \Sigma$
of codimension $q+1$, each $\mathfrak{u}\in \sigma^{\perp} \cap M$ gives a rational function $\chi^\mathfrak{u}$ on $V(\sigma)$.  The subspace $R^T_{q+1}(X_{\Sigma})$ is generated by
$\mathrm{div}\chi^\mathfrak{u}$ except the coefficients are taken mod $2$.
  That is, the coefficient of $[V(\tau)]$ in
the relation $\mathrm{div}\chi^\mathfrak{u}$ is
\begin{equation}\label{mod2}
 <\mathfrak{u} ,\mathfrak{n}_{\sigma,\tau}> \, \mathrm{mod}\,2
 \end{equation} where $\mathfrak{n}_{\sigma,\tau}$ is a lattice generator of $\nf{ \tau
\cap N}{ \sigma \cap N}$.
\begin{proposition}
For each $q$ with $0\leq q \leq d$ we have $$H_{qq}(\Sigma)\cong A^T_q(X_{\Sigma}),$$ where $A^T_q(X_{\Sigma})$ is the $q$th $\mathbb{Z}_2$ torus invariant Chow group of $X_{\Sigma}$.
\end{proposition}
To prove the proposition, let $\sigma < \tau$ in $\Sigma$ with $\mathrm{codim}\, \sigma=q+1$ and $\mathrm{codim}\, \tau=q$ and we consider the restriction map
\begin{equation}\label{diffj}
\rho_{\tau,\sigma} :\wedge^q N(\sigma)\longrightarrow \wedge^q N(\tau)
\end{equation}
for the cosheaf $\wedge^q\mathscr{E}$.
Choose a basis $\mathcal{B}:=\{t_1,t_2, \cdots , t_q\}$ for $N(\tau)$ and extend $\mathcal{B}$ to a basis
$\mathcal{B}^{\prime} :=\mathcal{B} \cup \{n_{\sigma,\tau}\}$ for $N(\sigma)$, where $n_{\sigma,\tau}$ is the nonzero element in $\nf{N_{\tau}}{N_{\sigma}}$.
We determine the map $\rho_{\tau,\sigma}$ in (\ref{diffj}) using the basis elements of $\wedge^qN(\sigma)$:  \begin{eqnarray*}
&\rho_{\tau,\sigma}(t_1\wedge t_2 \wedge \cdots \wedge t_q)&=t_1\wedge t_2 \wedge \cdots \wedge t_q \\
&\rho_{\tau,\sigma} (n_{\sigma,\tau}\wedge * \wedge \cdots \wedge *)&=0
\end{eqnarray*} where $t_1\wedge t_2 \wedge \cdots \wedge t_q$ is the generator of $\wedge^qN(\tau)$ and the $*$-entries can be any elements of $\mathcal{B}^{\prime}$.
Since $N(\sigma)$ has rank $q+1$, we have an isomorphism $\wedge^qN(\sigma) \cong N(\sigma)$ which identifies a $q$-tuple of basis elements
$b_1\wedge b_2 \wedge \cdots \wedge b_q$ to the one element of $\mathcal{B}^{\prime}\backslash \{b_1,b_2, \cdots ,b_q\}$.
Hence we can interpret the map $\rho_{\tau,\sigma}$ in (\ref{diffj}) as a map $N(\sigma)\longrightarrow \mathbb{Z}_2[\tau]$ which sends
the nonzero element $n_{\sigma,\tau}\in \nf{N_{\tau}}{N_{\sigma}}$  to $1\in \mathbb{Z}_2[\tau]$ and any $t_k$ to zero. \\
\indent Next, we consider the relations in  the $\mathbb{Z}_2$ torus invariant Chow groups $A_q^{T}(X_{\Sigma})$ which are given by $\mathrm{div}\chi^{\mathfrak{u}}$ as discussed above.  The element of $\mathbb{Z}_2$ in (\ref{mod2}) is the same as the element when we take $u$, the image of $\mathfrak{u}$ in $M(\sigma)=\nf{\sigma^{\perp} \cap M}{ \sigma^{\perp} \cap 2M}$ and $n_{\sigma,\tau}$, the image of $\mathfrak{n}_{\sigma,\tau}$ in $ \nf{N_{\tau}}{N_{\sigma}}$.
  Next, we have $M(\sigma) \cong \mathrm{Hom}(N(\sigma), \mathbb{Z}_2)$ and hence the subspace $R^T_{q+1}(X_{\Sigma})$ is generated by the
image of
\begin{eqnarray}\label{diffchow}
&\mathrm{Hom}(N(\sigma), \mathbb{Z}_2) \longrightarrow \mathbb{Z}_2[\tau]
\\& u \longmapsto   < u,n_{\sigma,\tau}> \mathrm{mod} \, 2.
  \end{eqnarray}
 We take the dual basis for $\mathcal{B}^{\prime}$ as a basis for $M(\sigma)$ which consists of elements $\{t^*_1,t^*_2, \cdots , t^*_q, n^*_{\sigma,\tau} \}$.
We determine the map (\ref{diffchow}) using this basis, and we see the image of (\ref{diffchow}) in $\mathbb{Z}_2[\tau]$ is equal to
the image of the differential (\ref{diffj}) in $\wedge^qN(\tau)\cong  \mathbb{Z}_2[\tau]$.  Moreover since  $C_{q-1}(\wedge^{q}\mathscr{E})=0$, we have  $$\mathrm{ker}(C_q(\wedge^q\mathscr{E})\longrightarrow C_{q-1}(\wedge^{q}\mathscr{E})) = C_q(\wedge^q\mathscr{E})$$
and as $C_q(\wedge^q\mathscr{E})\cong Z^T_q(X_{\Sigma})$ we obtain
\begin{eqnarray*}
\frac{C_q(\wedge^q\mathscr{E})}{\mathrm{im}(C_{q+1}(\wedge^q\mathscr{E})\longrightarrow C_q(\wedge^q\mathscr{E}))}&& \cong \quad \frac{Z^T_q(X_{\Sigma})}{R^T_{q+1}(X_{\Sigma})}\\
H_{qq}(\Sigma)&& \cong \quad A_q^{T}(X_{\Sigma}),
\end{eqnarray*} which proves the proposition.

Next, we recall the natural map
$$
f_q^{\mathbb{R}} :A^T_{q}(X_{\Sigma})\cong H_{qq}(\Sigma)\longrightarrow H_{q}(X_{\Sigma}(\mathbb{R}))
$$  given in (\ref{fmap}) which arises from the edge homomorphisms of the spectral sequence $(\widetilde{E}^r, \widetilde{d}^r)$.
  Thinking cellularly, an algebraic cycle $[V(\beta)]$ in $Z^T_q(X_{\Sigma})$ is represented by the sum of all $2^q$ cells in the orbit $O_{\beta}(\mathbb{R})$.
Hence, if  $g$ is a $q$ dimensional face of $\Delta$ with $\beta$ the cone dual to $g$ then
$f_q^{\mathbb{R}}\left( [V(\beta)] \right)$ is the homology class of $C_{\beta}$ where $C_{\beta}:= \sum_{t\in N(\beta)}(g,t)$ is the cellular chain in $C_{q}(X_{\Sigma})$ obtained by adding all $2^q$ copies of $g$ in $X_{\Sigma}(\mathbb{R})$.

 We say a cone $\sigma \subset N$ is \emph{$\mathbb{Z}_2$ regular} provided the image  in $\nf{N}{2N}$ of the rays of $\sigma$ forms a basis for the $\mathbb{Z}_2$ vector space $N_{\sigma}$.  Bihan et al. [Bih] show that if $\Sigma$ consists of $\mathbb{Z}_2$ regular cones then the map $f_q^{\mathbb{R}}$ is an isomorphism for $0\leq q \leq d$.  The counterexample in [How] shows that in general $f_q^{\mathbb{R}}$ is  neither injective nor surjective.  Next, we determine $A_q^{T}(X_{\Sigma})$
in some cases.  The following lemma will be useful in Section \ref{chowsec}.
\begin{lemma}\label{chowlemma}
Assume for each $\tau \in \Sigma(d-k-1)$ we have $V(\tau)(\mathbb{R})$ is $T$-homeomorphic to $\mathbb{RP}^{k+1}$.  If $ q\leq k$  then $A_q^{T}(X_{\Sigma})$ is generated by the orbit closure of any $q$ dimensional torus orbit in $X_{\Sigma}$.
\end{lemma}
To prove the lemma, we first look at the relations coming from a single cone $\sigma$ of codimension $q+1$.  Since $\sigma$ is a face of a cone in $\Sigma(d-k-1)$, $V(\sigma)(\mathbb{R})$ is $T$-homeomorphic to $\mathbb{RP}^{q+1}$.  Thus, $M(\sigma)\cong \mathrm{Hom}(N(\sigma),\mathbb{Z}_2)$ is generated by elements $a_1, a_2, \cdots ,a_{q+1}$ and $\sigma$ is contained in the $q+2$ cones $\omega_1, \omega_2 \cdots , \omega_{q+2}$ of codimension $q$, where
$$ n_{\sigma, \omega_i} =\left\{ \begin{array}{ll}a_i & \mathrm{if} \quad 1\leq i \leq q+1 \\a_1+a_2+ \cdots +a_{q+1} & \mathrm{if} \quad i=q+2.  \end{array} \right. $$
Hence, the map
\begin{equation}\label{omegamap}
M(\sigma) \longrightarrow Z^T_q(X_{\Sigma})
\end{equation} is given by basis elements as follows.
\begin{eqnarray*}
&a_1\longmapsto [V(\omega_1)]+[V(\omega_{q+2})]\\
&a_2\longmapsto [V(\omega_2)]+[V(\omega_{q+2})]\\
&\cdots \\
&a_{q+1}\longmapsto [V(\omega_{q+1})]+[V(\omega_{q+2})]
\end{eqnarray*}
 Extending linearly to all of $M(\sigma)$, we see that in the image of (\ref{omegamap}) we obtain the sum of any even number of $[V(\omega_i)]$'s.  For any $r$ with $1\leq r \leq q+2$, the cycle $[V(\omega_r)]$ is not in the image of (\ref{omegamap}).  Moreover, if $C$ is a sum of an odd number of $[V(\omega_i)]$'s then $C+[V(\omega_r)]$ is in the image of (\ref{omegamap}). Hence for any $r$, the cycle $[V(\omega_r)]$ is the generator of $\displaystyle{\mathrm{coker}(M(\sigma)\longrightarrow \bigoplus_{\sigma < \omega_i}[V(\omega_i)])}.$
 Next, we consider the map
 \begin{equation} \label{omegamap2}
 \bigoplus_{\sigma \in \Sigma(d-q-1)}M(\sigma) \longrightarrow \bigoplus_{\omega \in \Sigma(d-q)}[V(\omega)] , \end{equation}
 where the map on each $M(\sigma)$ is described above.
 \begin{claim}
 If $\omega$ and $\omega^{\prime}$ are codimension $q$ cones in $\Sigma$ then $[V(\omega)]+[V(\omega^{\prime})]$ is in the image of (\ref{omegamap2}).
 \end{claim}
  To see the claim, note that if $\omega \cap \omega^{\prime}$ is a codimension $q+1$ cone then $[V(\omega)]+[V(\omega^{\prime})]$ is in the image of (\ref{omegamap2}), as discussed above.  Else, find a sequence of codimension $q$ and codimension $ q+1$ cones
  $$
  \xymatrix{\omega=\omega_0 \ar@{:}[dr]&& \omega_1 \ar@{:}[dl] \ar@{:}[dr]&& \omega_2 \ar@{:}[dl] &\cdots &\omega_{n-1} \ar@{:}[dr] &&\omega_n=\omega^{\prime} \ar@{:}[dl]\\ & \sigma_{01} && \sigma_{12} &&\cdots && \sigma_{n-1n}&
  }$$ where $\sigma_{i-1i}< \omega_{i-1}$ and $\sigma_{i-1i}< \omega_{i}$ for $1\leq i \leq n$.
  There exists $u_{i-1i}\in M(\sigma_{i-1i})$ with $\mathrm{div}(\chi^{\mathfrak{u}_{i-1i}})= [V(\omega_{i-1})]+[V(\omega_{i})]$.  Adding, we obtain $$\sum_{i=1}^{n}u_{i-1i} \in \bigoplus_{\sigma \in \Sigma(d-q-1)}M(\sigma)$$ and $$\sum_{i=1}^{n}u_{i-1i} \longmapsto [V(\omega)]+[V(\omega^{\prime})]$$
  which proves the claim.  Thus, the sum of any even number of $[V(\omega)]$'s is in the image  of (\ref{omegamap2}) and the cycle $[V(\omega^{\prime})]$ is not.  Moreover, if $C$ is a chain in $Z^T_q(X_{\Sigma})$ with an odd number of $[V(\omega)]$'s then $C+[V(\omega^{\prime})]$ is in the image of (\ref{omegamap2}).  Hence, $A_q^T(X_{\Sigma})$ is generated by $[V(\omega^{\prime})]$ where $\omega^{\prime}$ is any cone of codimension $q$.  This proves Lemma \ref{chowlemma}.
\subsection{The right-most column $H_{dq}(\Sigma)$}\label{rightsec}
Let $\{r_1, r_2, \cdots , r_k \}$ be the rays of $\Sigma$ and $R_i$ the image of $r_i$ in $\nf{N}{2N}\cong (\mathbb{Z}_2)^d$. 
\begin{proposition}
 Suppose
the rank of the $\mathbb{Z}_2$ vector space $V:=\mathrm{span}(R_1, R_2, \cdots , R_k)$ is $s\leq d$. 
Then, $$\mathrm{rank}H_{dq}(\Sigma)={d-s \choose q-s} \quad \mbox{and} \quad \overline{d}^r_{d,q}=0\quad \mbox{for} \quad r\ge 1.$$
\end{proposition}
 To prove the proposition, we choose a basis for $V$ consisting of a subset of the $R_i$, and reorder if needed so that $\{R_1, R_2, \cdots , R_s\}$ is a basis for $V$.  Note that the elements of the form $$R_i\wedge*\wedge * \wedge \cdots \wedge *$$
generate $$\mathrm{ker}(  \wedge^q\nf{N}{2N} \longrightarrow \wedge^q N(r_i)),$$ where the $*$-entries are any elements of $\nf{N}{2N}$.
Hence,  $$\mathrm{ker}(\wedge^q\nf{N}{2N} \longrightarrow \bigoplus_{i=1}^s \wedge^q N(r_i))$$ is generated by elements of the following form.
$$R_1\wedge R_2\wedge \cdots \wedge R_s \wedge *\wedge * \wedge \cdots \wedge *$$  This is a subspace of $\wedge^q\nf{N}{2N}$ of rank
${d-s \choose q-s}$, where  ${n \choose k}=0$ if $k<0$.  Hence, the kernel of the boundary map $$\overline{d}^0_{d,q}: \wedge^q \nf{N}{2N} \longrightarrow \bigoplus_{r\in \Sigma(1)}\wedge^q N(r)$$
has rank at most ${d-s \choose q-s}$.  This shows
\begin{equation}\label{leq}
b_d(X_{\Sigma}(\mathbb{R})) \leq \sum_{q=s}^d {d-s \choose q-s}=2^{d-s}
\end{equation}
 and equality holds if and only if for each $q$, $\mathrm{ker}\overline{d}^0_{d,q}$ has rank ${d-s \choose q-s}$ and all higher differentials $\overline{d}^r_{d,q}$, $r\ge 1$ with source the rightmost column are zero.   Next, we note that after identifying all facets of the $2^d$ copies of $\Delta$, we are left with $2^{d-s}$ components which shows
\begin{equation}\label{ge}
b_d(X_{\Sigma}(\mathbb{R}))\ge 2^{d-s}.
\end{equation}
Combining (\ref{leq}) and (\ref{ge}) we have equality yielding $b_d(X_{\Sigma}(\mathbb{R}))= 2^{d-s}, \mathrm{rank}H_{dq}(\Sigma)={d-s \choose q-s}$, and the higher boundaries $\overline{d}^r_{d,q}$, $r\ge 1$ are zero.

\section{$\mathbb{Z}_2$ Hodge spaces for reflexive polytopes}\label{reflexsection}
\subsection{A correspondence between sheaves and cosheaves} \label{corrsec}
A \emph{reflexive} polytope is a lattice polytope $\Delta$ with $0 \in \mathrm{int} \Delta$ and such that the
polar polytope $\Delta^*$ is also a lattice polytope. A discussion of reflexive polytopes can be found in [Bat1 \textsection 4.1].
Throughout this section we will use the following notation.
\begin{eqnarray*}  \Delta \subset M &&  \mbox{a reflexive polytope}\\ \Delta^*\subset N && \mbox{the polar polytope of } \Delta
\\ \Sigma\subset N &&\mbox{the normal fan of } \Delta \\ && (=\mbox{the face fan of } \Delta^*) \\
 \Sigma^* \subset M && \mbox{the normal fan of } \Delta^* \\&& (=\mbox{the face fan of } \Delta)
\end{eqnarray*}
If $\tau\in \Sigma$ and $\mathrm{dim}\, \tau=j>0$ then $\tau=\mathrm{poshull}f^*$ where $f^*<\Delta^*$ is a face of dimension $j-1$.
We define $\tau^*$ to be the cone in $\Sigma^*$ of dimension $d-j+1$ which is dual to $f$.  The correspondence
$$ \tau\in  \Sigma \longleftrightarrow \tau^* \in \Sigma^*$$ gives a one-to-one inclusion reversing correspondence between the
cones in $\Sigma$ of positive dimension
and the positive dimensional cones in $\Sigma^*.$ We show this correspondence by dimension below.
 \begin{center}
 $\xymatrix{
 \mathrm{Cone}\, \mathrm{in} \, \Sigma: & 1 \ar@{<->}[d]& 2 \ar@{<->}[d] &3  \ar@{<->}[d]& \cdots & d-2  \ar@{<->}[d]& d-1  \ar@{<->}[d]& d  \ar@{<->}[d] \\ \mathrm{Cone} \, \mathrm{in} \, \Sigma^*:
& d & d-1 & d-2 & \cdots & 3 & 2 & 1
 \\ } $
 \end{center}
Let $\mathscr{H}$ be a sheaf on the fan $\Sigma^*$.  We use the correspondence above to create a cosheaf $\widehat{\mathscr{H}}$ on
$\Sigma$ by defining for $\tau \in \Sigma$
\begin{equation}\label{stalk}
\widehat{H}_{\tau}:=\left\{ \begin{array}{ll}H_{\tau^*}&\mathrm{if} \, \mathrm{dim}\, \tau>0\\
 0 &\mathrm{if} \, \mathrm{dim}\, \tau=0 \end{array} \right.
 \end{equation}
with restriction map for $\sigma <\tau $
 \begin{equation}\label{resmap}
 \widehat{\rho}_{\tau,  \sigma}:=\left\{ \begin{array}{ll}\rho_{\tau^* ,\sigma^*}&\mathrm{if} \, \mathrm{dim} \, \sigma >0\\
 0 &\mathrm{if} \, \mathrm{dim} \, \sigma=0 \end{array} \right.
 \end{equation}
  where $\rho_{\tau^* ,\sigma^*}$ is a restriction map for the sheaf $\mathscr{H}$.
 Using (\ref{stalk}) we have equality of chain groups for the sheaf $\mathscr{H}$ on $\Sigma^*$ and the cosheaf $\widehat{\mathscr{H}}$ on $\Sigma$, as depicted below.
 $$\xymatrix{ C_{d-1}(\widehat{\mathscr{H}}\, )\ar[r] \ar@{=}[d] & C_{d-2}(\widehat{\mathscr{H}}\, )  \ar[r] \ar@{=}[d] & \cdots \ar[r] & C_{2}(\widehat{\mathscr{H}}\, )\ar[r] \ar@{=}[d] & C_{1}(\widehat{\mathscr{H}}\, )\ar[r] \ar@{=}[d] & C_{0}(\widehat{\mathscr{H}}\,) \ar@{=}[d] \\
 C^0(\mathscr{H}) \ar[r] & C^1(\mathscr{H})\ar[r] & \cdots \ar[r] &  C^{d-3}(\mathscr{H}) \ar[r] & C^{d-2}(\mathscr{H})  \ar[r] & C^{d-1}(\mathscr{H}) \\
  }$$ From (\ref{resmap}) the horizontal sheaf and cosheaf boundary maps are equal.  Hence we may equate sheaf cohomology groups for $\mathscr{H}$ on $\Sigma^*$ with cosheaf homology groups for $\widehat{\mathscr{H}}$ on $\Sigma$.
  $$ H^p(\mathscr{H}) \cong H_{d-p-1}(\widehat{\mathscr{H}}\, ) \quad 1\leq p \leq d-2 $$
  One more construction which we will use in the sequel is that of the cosheaf $\mathscr{A}^{\circ}$.  If $\mathscr{A}$ is a cosheaf on $\Sigma$ then 
$\mathscr{A}^{\circ}$ is defined by
  $$A_{\sigma}^{\circ}:=\left\{ \begin{array}{ll} A_{\sigma} &\mbox{if} \quad \mathrm{dim} \, \sigma>0\\ 0 &\mbox{if} \quad \mathrm{dim} \,\sigma=0\\ \end{array}  \right. $$ 
with the restriction map for $\sigma < \tau$ defined by
    $$\rho_{\tau,\sigma}^{\circ}:=\left\{ \begin{array}{ll} \rho_{\tau,\sigma}&\mbox{if} \quad \mathrm{dim} \, \sigma>0\\ 0 &\mbox{if} \quad \mathrm{dim} \,\sigma=0. \\ \end{array}  \right.  $$
    Note that by definition $$C_p(\mathscr{A}) =C_p( \mathscr{A}^{\circ}) \quad \mbox{for} \quad p\leq d-1,$$
    and hence we have $$H_p(\mathscr{A})=H_p(\mathscr{A}^{\circ}) \quad \mbox{for}\quad p\leq d-2 .$$
\subsection{The cosheaves $\widehat{\mathscr{F}}$, $\widehat{\mathscr{G}}$, and $\mathscr{C}$ on $\Sigma$}\label{fcsection}
The sheaf $\mathscr{F}$ on $\Sigma^*$ is defined as follows.
For $\sigma^* \in \Sigma^*$ the stalk is
$$F_{\sigma^*}:= \frac{(\sigma^*)^{\perp} \cap N}{(\sigma^*)^{\perp} \cap 2N}$$ and the face restriction map
$$\xymatrix{ \rho_{\sigma^*, \tau^*}: F_{\sigma^*} \ar[r]^{\quad \subset}& F_{\tau^*}}$$ is given by inclusion for $\tau^*< \sigma^*$ in $\Sigma^*$.  The sheaf $\mathscr{G}$ on $\Sigma^*$ is then defined to be the cokernel of the inclusion $\mathscr{F}\hookrightarrow \nf{N}{2N}$ so that
\begin{equation}\label{brion}
0\longrightarrow \mathscr{F} \longrightarrow \nf{N}{2N} \longrightarrow \mathscr{G} \longrightarrow 0
\end{equation} is a short exact sequence of sheaves on $\Sigma^*$.
\begin{claim} For $\tau \in \Sigma, \tau^* \in \Sigma^*$ as in the previous section with $\mathrm{dim}\, \tau=l >0$,
we have the following containment of $\mathbb{Z}_2$ vector spaces
\begin{equation} \label{shorter}\frac{(\tau^*)^{\perp} \cap N}{(\tau^*)^{\perp} \cap 2N} \subset \frac{\tau \cap N}{\tau \cap 2N} .\end{equation}
\end{claim}
Note that the $\mathbb{Z}_2$ vector space on the left has rank $l-1$ while the one on the right has rank $l$.
To prove the claim, let $f<\Delta$ and $f^* < \Delta^*$ be such that $$ \begin{array}{rlrl}\tau &= \mathrm{poshull}f^* &\quad \tau^* &= \mathrm{poshull}f
\\ & = \mathrm{poshull}\{q_1,q_2, \cdots ,q_s\} \quad & & =\mathrm{poshull}\{p_1,p_2, \cdots ,p_k\} .\end{array}$$
Then we have $<p_i, q_j>=-1 \quad \forall \, i,j.$
This gives $(\tau^*)^{\perp}=\mathrm{span}\{q_i-q_j \, | \, 1\leq i,j\leq s\}$ and the claim follows.  Next, we note that the inclusion (\ref{shorter}) is compatible with the face restriction maps for cones in $\Sigma$.  Thus, we obtain an injective homomorphism of cosheaves $\widehat{\mathscr{F}} \hookrightarrow \mathscr{N}. $  We define the cosheaf $\mathscr{C}$ to be the cokernel of this homomorphism so that
 \begin{equation} \label{fcseq} 0\longrightarrow \widehat{\mathscr{F}} \longrightarrow \mathscr{N}  \longrightarrow \mathscr{C} \longrightarrow 0 \end{equation}
is an exact sequence of cosheaves on $\Sigma$.  Note that if $\mathrm{dim} \,\sigma=0$ then the stalks $\widehat{F}_{\sigma}$, $N_{\sigma}$, and $C_{\sigma}$ are all zero.   Moreover, for each $\sigma \in \Sigma$ of positive dimension the stalk $C_{\sigma}$ is a rank one $\mathbb{Z}_2$ vector space.
\begin{lemma} \label{fclemma}
For $\sigma < \tau$ in $\Sigma$ with $\mathrm{dim} \,\sigma>0$ the map $C_{\sigma} \longrightarrow C_{\tau}$ is the identity.
\end{lemma}
We prove the lemma by contradiction.  Assume we have $\tau_1< \tau_2$ in $\Sigma$ and $C_{\tau_1} \longrightarrow C_{\tau_2}$ is the zero map.  We have the follow diagram
$$ \xymatrix{ 0 \ar[r] & \frac{(\tau_1^*)^{\perp} \cap N}{(\tau_1^*)^{\perp} \cap 2N} \ar[d] \ar[r] & \frac{\tau_1 \cap N}{\tau_1 \cap 2N} \ar[d] \ar[r] & C_{\tau_1} \ar[d]^{\phi} \ar[r] & 0 \\
0 \ar[r] &\frac{(\tau_2^*)^{\perp} \cap N}{(\tau_2^*)^{\perp} \cap 2N} \ar[r] &\frac{\tau_2 \cap N}{\tau_2 \cap 2N}  \ar[r] & C_{\tau_2} \ar[r] & 0
}$$ where the rows are short exact sequences and vertical maps are cosheaf restriction maps.  We have assumed $\phi=0$ and hence
\begin{equation}\label{incl}
\frac{\tau_1 \cap N}{\tau_1 \cap 2N} \subset \frac{(\tau_2^*)^{\perp} \cap N}{(\tau_2^*)^{\perp} \cap 2N}.
\end{equation}
Let $r_1$ be the first lattice point on a ray of $\tau_1 \subset \tau_2$ and $r_2$ be the first lattice point on a ray of $\tau_2^*$.  Due to the inclusion (\ref{incl}) we have $<r_1,r_2>= 0\, (\mathrm{mod} \, 2)$.  This is a contradiction because $<r_1,r_2>=-1$.  Thus we have Lemma \ref{fclemma}. \\
\indent Next, note that 
\begin{equation} \label{yousuck}0\longrightarrow \widehat{\mathscr{F}}\longrightarrow \widehat{\nf{N}{2N}}\longrightarrow \widehat{\mathscr{G}} \longrightarrow 0 \end{equation}
is a short exact sequence of cosheaves of $\Sigma$ where $\widehat{\nf{N}{2N}}=\nf{N}{2N}^{\circ}$ as cosheaves on $\Sigma$.  We can combine (\ref{yousuck}) with (\ref{fcseq}) into the following commutative diagram \vspace{-.45in}
 \begin{equation*}\xymatrix{ \\ & &  & 0\ar[d] \\ &0 \ar[d] &0 \ar[d]&\mathscr{K} \ar[d] & \\ 0 \ar[r] &
 \widehat{\mathscr{F}}\ar[d] \ar[r] & \widehat{\nf{N}{2N}}\ar[d] \ar[r] & \widehat{\mathscr{G}} \ar[d]^{\Phi} \ar[r]& 0 \\
 0 \ar[r] & \mathscr{N}\ar[d] \ar[r] &\nf{N}{2N}^{\circ} \ar[d] \ar[r] & \mathscr{E}^{\circ} \ar[d] \ar[r]& 0\\
 & \mathscr{C} \ar[d] & 0 & 0\\ & 0 & &
}\end{equation*} 
where the rows and columns are exact and $\mathscr{K}:=\mathrm{ker}\Phi$. \\
\indent Next, we use the Snake Lemma on the commutative diagram to obtain $\mathscr{K}\cong \mathscr{C}$ as cosheaves yielding
 \begin{equation} \label{circexact}0\longrightarrow \mathscr{C} \longrightarrow \widehat{\mathscr{G}} \longrightarrow \mathscr{E}^{\circ} \longrightarrow 0
\end{equation} is a short exact sequence of cosheaves on $\Sigma$.   As $(\wedge^q\mathscr{E})^{\circ}= \wedge^q\mathscr{E}^{\circ}$ we have  $H_p(\wedge^q \mathscr{E}) =H_p(\wedge^q \mathscr{E}^{\circ})$ for $p\leq d-2.$
\subsection{Vanishing of  the homology groups $H_p(\wedge^k\widehat{\mathscr{G}}\, )$}\label{gsec}
We remind the reader that a cone
$\sigma \in \Sigma$ is \emph{$\mathbb{Z}_2$ regular} provided the image in $\nf{N}{2N}$ of the rays of $\sigma$ forms a basis for $N_{\sigma}$. 
\begin{lemma} \label{gvanish}
Assume the cones in $\Sigma^*$ of dimension at most $e$ are $\mathbb{Z}_2$ regular.  Then $$H_p(\wedge^k \widehat{\mathscr{G}}\, )=0 \quad \mbox{for} \quad 1\leq p < e-1.$$
\end{lemma}
First, we will assume that $e=d$.  Following [Bri \textsection1.2], as the cones in $\Sigma^*$ are $\mathbb{Z}_2$ regular the sheaf $\mathscr{G}$ on $\Sigma^*$ can be written as follows.
$$\mathscr{G}=\bigoplus_{r_i\in \Sigma^*(1)}\mathscr{G}(r_i)$$ is a direct sum of the sheaves $\mathscr{G}(r_i)$ with
$$G(r_i)_{\tau^*}= \left\{ \begin{array}{ll}\mathbb{Z}_2 & \mathrm{if}\, r_i \in \tau^*(1)\\
0 &\mathrm{else}.
\end{array} \right.$$
Moreover, we have $$\wedge^k\mathscr{G}= \bigoplus_{r_1,r_2, \cdots ,r_k \, \mathrm{distinct}}\mathscr{G}(r_1, r_2, \cdots , r_k)$$ is the direct sum of sheaves $\mathscr{G}(r_1, r_2, \cdots , r_k)$ with
$$G(r_1, r_2, \cdots , r_k)_{\tau^*}= \left\{ \begin{array}{ll}\mathbb{Z}_2 & \mathrm{if}\, \, r_1, r_2, \cdots ,r_k \in \tau^*(1)\\
0 &\mathrm{else}.
\end{array} \right.$$
The proof of the proposition in Section 1.2 of [Bri] yields
$H^p(\mathscr{G}(r_1, r_2, \cdots ,r_k))=0$ for  $p>0$ because $\Sigma^*$ is a complete fan.
Hence we have
\begin{eqnarray*}
H^p(\wedge^k\mathscr{G})& &= H^p(\bigoplus_{r_1,r_2, \cdots ,r_k \, \mathrm{distinct}}\mathscr{G}(r_1, r_2, \cdots , r_k))\\
 && \cong \bigoplus_{r_1,r_2, \cdots ,r_k \, \mathrm{distinct}}H^p(\mathscr{G}(r_1, r_2, \cdots , r_k)) \\
 && =0 \quad \mathrm{for} \quad p>0. \end{eqnarray*}
 Hence, as $\wedge^k \widehat{ \mathscr{G}}= \widehat{\wedge^k \mathscr{G}}$ we have $H_p(\wedge^k\widehat{\mathscr{G}}\, )=H^{d-p-1}(\wedge^k \mathscr{G})=0 $ for $p$ such that $1\leq p \leq d-2$ and $d-p-1>0$, which proves Lemma \ref{gvanish} when $e=d$.

Next, assume $e<d$.  Let $\Sigma^*_{\leq e} := \bigcup_{i\leq e} \Sigma^*(i)$ be the subfan of $\Sigma^*$ consisting of the cones of dimension at most $e$, and let $\mathscr{G}^{\prime}$
be the restriction of $\mathscr{G}$ to $\Sigma^*_{\leq e}$.  Then as $\Sigma^*_{\leq e}$ consists of $\mathbb{Z}_2$ regular cones, we have
$$\mathscr{G}^{\prime}=\bigoplus_{r_i\in \Sigma^*(1)}\mathscr{G}(r_i)^{\prime}$$ where $\mathscr{G}(r_i)^{\prime}$ is the restriction of $\mathscr{G}(r_i)$ to $\Sigma^*_{\leq e}.$
Moreover, we have $$\wedge^k\mathscr{G}^{\prime}= \bigoplus_{r_1,r_2, \cdots ,r_k \, \mathrm{distinct}}\mathscr{G}(r_1, r_2, \cdots , r_k)^{\prime}.$$
By definition, for each $k$ we have the equality $H^p(\wedge^k\mathscr{G})=H^p(\wedge^k\mathscr{G}^{\prime})$ for  $p> d-e$ which gives $H^p(\wedge^k\mathscr{G})=0$ for $p>d-e$.
Thus, $H_p(\wedge^k\widehat{ \mathscr{G}}\, )=0$ for $p$ such that $1\leq p \leq d-2$ and $d-p-1>d-e$, and we have Lemma \ref{gvanish}.
\subsection{Vanishing of the $\mathbb{Z}_2$ Hodge spaces $H_{pq}(\Sigma)$} \label{vansec}
\begin{theorem}  Assume the cones in $\Sigma^*$ of dimension at most $e$ are $\mathbb{Z}_2$ regular.  Then $$H_p(\wedge^q\mathscr{E})=0\quad \mbox{for} \quad q < p < e-1.$$
\end{theorem}
To prove the theorem, we use the short exact sequence
 $0\longrightarrow \mathscr{C} \longrightarrow \widehat{\mathscr{G}} \longrightarrow \mathscr{E}^{\circ} \longrightarrow 0$ of cosheaves of $\Sigma$ and the associated degree $q$ sequence from (\ref{koscsh})
 $$ 0\longrightarrow S^q\mathscr{C}\longrightarrow \wedge^1\widehat{\mathscr{G}}\otimes S^{q-1}\mathscr{C}\longrightarrow \cdots \longrightarrow
\wedge^{q-1}\widehat{\mathscr{G}}\otimes S^1\mathscr{C}\longrightarrow \wedge^q\widehat{\mathscr{G}}\longrightarrow \wedge^q\mathscr{E}^{\circ} \longrightarrow 0 $$
which we break into short exact sequences below
\begin{eqnarray}
 & 0\longrightarrow S^q\mathscr{C}\longrightarrow \wedge^1\widehat{\mathscr{G}}\otimes S^{q-1}\mathscr{C} \longrightarrow W_{1} \longrightarrow 0 \label{eq1} \\
 & 0 \longrightarrow W_{1} \longrightarrow  \wedge^2\widehat{\mathscr{G}}\otimes S^{q-2}\mathscr{C} \longrightarrow W_{2} \longrightarrow 0 \label{eq2} \\
 & \cdots \\
 & 0 \longrightarrow W_{q-1} \longrightarrow  \wedge^{q}\widehat{\mathscr{G}} \longrightarrow \wedge^q\mathscr{E}^{\circ} \longrightarrow 0. \label{eqq}
\end{eqnarray}
These induce long exact sequences on homology groups.  From Lemma \ref{fclemma} we see that 
\begin{eqnarray*}\mathscr{C}\cong (\mathbb{Z}_2)^{\circ} \\ S^{q-k}\mathscr{C} \cong (\mathbb{Z}_2)^{\circ} \end{eqnarray*}  where $\mathbb{Z}_2$ is the constant cosheaf on $\Sigma$.  Thus
 we have $$H_p(\wedge^k\widehat{\mathscr{G}}\otimes S^{q-k}\mathscr{C}) \cong H_p(\wedge^k\widehat{\mathscr{G}}\, ) \otimes S^{q-k}\mathscr{C}=0 \quad \mbox{for} \quad 1\leq p <e-1,$$
  where we have used Lemma \ref{gvanish}.
We begin with the long exact sequence induced from (\ref{eq1}).
$$\cdots \longrightarrow H_p( \wedge^1\widehat{\mathscr{G}}\otimes S^{q-1}\mathscr{C} )\longrightarrow H_p(W_1) \longrightarrow H_{p-1}( S^q\mathscr{C}) \longrightarrow \cdots $$
We have
\begin{eqnarray*}H_p( \wedge^1\widehat{\mathscr{G}}\otimes S^{q-1}\mathscr{C} )=0 && \mathrm{for} \quad 1\leq p <e-1\\
H_{p-1}( S^q\mathscr{C})=0 && \mathrm{for} \quad 1\leq p-1 <e-1 \end{eqnarray*} and hence $H_p(W_1)=0$ for $1<p<e-1$.  Next we use the long exact sequence induced from (\ref{eq2}).
$$\cdots \longrightarrow H_p( \wedge^2\widehat{\mathscr{G}}\otimes S^{q-2}\mathscr{C} )\longrightarrow H_p(W_2) \longrightarrow H_{p-1}( W_1) \longrightarrow \cdots $$
We have
\begin{eqnarray*}
H_p( \wedge^2\widehat{\mathscr{G}}\otimes S^{q-2}\mathscr{C} )=0 && \mathrm{for} \quad 1\leq p <e-1 \\
H_{p-1}(W_1)=0 && \mathrm{for} \quad 1<p-1<e-1
\end{eqnarray*}
and hence $H_p(W_2)=0$ for $2<p<e-1$.  We continue this process to obtain $H_p(\wedge^q\mathscr{E}^{\circ})=0$ for  $q < p < e-1.$  Moreover, as $e-1 \leq d-1$ we have $H_p(\wedge^q\mathscr{E})=H_p(\wedge^q \mathscr{E}^{\circ})$ for $q<p<e-1$ and the theorem holds.

\subsection{The diagonal entries $H_{qq}(\Sigma)$}\label{chowsec}
In this section, we are under the assumption that the cones in $\Sigma^*_{\leq e}$ are $\mathbb{Z}_2$ regular.
Let $\tau^* \in \Sigma^*(e)$, $\tau^*=\mathrm{poshull}\, f$ where $f=\mathrm{conv}\{p_1,p_2, \cdots ,p_e\}$.
\begin{lemma}\label{sublemma}
The toric subvariety $Y(\mathbb{R})$ of $X_{\Sigma}(\mathbb{R})$ defined by the face $f$ of $\Delta$ is $T$-homeomorphic to $\mathbb{RP}^{e-1}$.
\end{lemma}
To prove the lemma, let $\Psi=\mathrm{conv}\{ 0, p_1, p_2, \cdots , p_e\} \subset M$.  We extend $\{p_1, p_2, \cdots ,p_e\}$ to a basis $\{p_1, p_2, \cdots ,p_e, t_1, t_2, \cdots , t_{d-e} \}$  for $M\otimes \mathbb{R}$ with  $\{t_1, t_2 ,\cdots ,t_{d-e}\}$ orthonormal. We have a map $M^{\prime} \longrightarrow M$ given by the $d \times d$ matrix $$A:= \left[\begin{matrix} | & |& & | & | & & | \\ p_1 & p_2& \cdots &  p_e & t_{1}&  \cdots & t_{d-e}\\  | & |& & | & | & & |  \end{matrix} \right]$$ which has odd determinant because $\tau^*$ is $\mathbb{Z}_2$ regular.  Moreover, we have $Ae_i =p_i$ for $1\leq i \leq e$ and $A$ takes the simplex $\Psi^{\prime}=\mathrm{conv}\{0, e_1, e_2, \cdots, e_e\}$ to the simplex $\Psi$.   The matrix $A^{*}$ gives a map $N \longrightarrow N^{\prime}$ which induces an isomorphism
\begin{equation}\label{nisom} \xymatrix{
\nf{N}{2N} \ar^{\cong}[r]& \nf{N^{\prime}}{2N^{\prime}}.  }\end{equation}
Let $r_i$ be the first integer point on the ray dual to the facet $\mathrm{conv}\{0, p_1, p_2, \cdots, \widehat{p_i}, \cdots , p_e\}$ of $\Psi$.  We have $0=<Ae_j, r_i>=<e_j, A^{*}r_i>$ for $j \in \{1,2, \cdots , \widehat{i}, \cdots e \}$, and hence the vector $A^{*}r_i$ lies on the ray dual to the facet $\mathrm{conv}\{0, e_1, e_2, \cdots, \widehat{e_i}, \cdots , e_e\}$ of $\Psi^{\prime}$.  This gives $A^{*}r_i=ke_i,$ where $k$ is an odd integer. (if $k$ were even, then the image of $r_i$ in $\nf{N}{2N}$ would map to $0 \in \nf{N^{\prime}}{2N^{\prime}}$ contradicting the isomorphism (\ref{nisom}) )
Similarly, if $r$ is the first integer point along the ray in $N$ dual to the facet
  $\mathrm{conv}\{ p_1, p_2, \cdots, p_e\}$ of $\Psi$ then $A^{*}r$ is an odd multiple of the vector $-e_1-e_2- \cdots - e_e$.
 
 Thus, the isomorphism in (\ref{nisom}) gives an isomorphism $\xymatrix{ N_{\sigma}\ar^{\cong}[r] &N_{\sigma^{\prime}} }$
where $\sigma$ is a cone in the normal fan of $\Psi$ and $\sigma^{\prime}$ is in the normal fan of $\Psi^{\prime}$.
 We obtain a $T$-homeomorphism between the real toric variety defined by $\Psi$ and $\mathbb{RP}^{e}$.  As $f<\Psi$ the toric variety $Y$ is  $T$-homeomorphic to $\mathbb{RP}^{e-1}$ and we obtain Lemma \ref{sublemma}.
Next, we use Lemma \ref{chowlemma} to arrive at the following proposition.
\begin{proposition}
If the cones in $\Sigma^*_{\leq e}$ are $\mathbb{Z}_2$ regular then for $q< e-1$ we have $H_{qq}(\Sigma)\, \, (\cong A^T_q(X_{\Sigma}))$ has rank $1$ and is
generated by the orbit closure of any $q$ dimensional torus orbit.
\end{proposition}
\section{Collapsing of the spectral sequence $\overline{E}^r$}\label{collsec}
In this section, we show that the spectral sequence $\overline{E}^r$ for $X_{\Sigma}$ collapses at $\overline{E}^1$ when the cones in $\Sigma^*$
are $\mathbb{Z}_2$ regular.  Combining the work in Section \ref{vansec} and Section \ref{chowsec}, we have that the
ranks of the entries in the $\overline{E}^1$ term are as follows
\begin{equation}\label{e1smooth}
q \begin{matrix} & & & & &  & & 1\\ &&&&& &*&*\\  &&&&&1&*&*\\ &&&&1&0&*&*\\  &&&1&0&0&*&* \\&&1&0&0&0&*&* \\
&1&0&0&0&0&*&*\\  1 & 0 & 0 & 0 & 0 & 0 & 0 & 0  \\ && & p &&& \end{matrix}
\end{equation}
where the $*$ entries are possibly nonzero and occur in the columns $p> d-2$.
\begin{remark} We can completely determine the $*$ entries.  From Section \ref{rightsec}, if $s$ is the rank of the image of the
rays of $\Sigma$ in $\nf{N}{2N}$ then the rank of $\overline{E}^1_{d,q}$ is the binomial coefficient
${d-s \choose q-s}$.  The Euler characteristic of the $q$th row is
$(-1)^qh_q$ where $h=(h_0,h_1,h_2, \cdots ,h_d)$ is the $h$-vector of the polytope $\Delta$.  Thus, we may determine the ranks of
the vector spaces $\overline{E}^1_{d-1,q}$, giving us knowledge of the ranks of all the entries in $\overline{E}^1_{p,q}$.
\end{remark}
By looking at the $\overline{E}^1$ term for $X_{\Sigma}$, we see that the only possible nonzero higher differentials have target
$\overline{E}^1_{d-2,d-2}$.  To show that the spectral sequence $\overline{E}^r$ for $X_{\Sigma}$ collapses at
$\overline{E}^1$, we need only show the following lemma.
\begin{lemma} \label{collapse}
The map $$f_{d-2}^{\mathbb{R}}:\overline{E}^1_{d-2,d-2} \longrightarrow H_{d-2}(X_{\Sigma}(\mathbb{R}))$$
from (\ref{fmap}) is nonzero.
\end{lemma}
Let $g$ be a $d-2$ dimensional face of $\Delta$ with $\beta= \mathrm{poshull}\, \{q_1, q_2, \cdots ,q_k \}$ the cone dual to $g$.
As mentioned in Section \ref{diagsec}, $f_{d-2}^{\mathbb{R}}$ sends $[V(\beta)]$ to the homology class of
 $C_{\beta}:= \sum_{t\in N(\beta)}(g,t)$.
 Suppose $C_{\beta}=\partial C$, $C =(g_1, t_1)+(g_2,t_2)+ \cdots +(g_r,t_r) \in C_{d-1}(X_{\Sigma})$ where $g_i <\Delta$,
$\sigma_i \in \Sigma$ is dual to $g_i$, and $t_i \in N(\sigma_i)$.  We include $N \subset \widetilde{N}$ where
$\widetilde{N} \cong \mathbb{Z}^{d+1}$ and $\Delta^*\times[-1,1] \subset \widetilde{N}$.  Note that $\Delta^*\times[-1,1] $ is
reflexive and the normal fan of $\Delta^*\times[-1,1]$ consists of $\mathbb{Z}_2$ regular cones.  Let $\Xi$ be the face fan
of $\Delta^*\times[-1,1]$.  Then, $\Xi$ is the normal fan of $B_{\Delta}:= (\Delta^*\times [-1,1])^*$ the bipyramid with base $\Delta$.
 We define the cone $\widetilde{\beta}\in \Xi$ to be the positive hull of the rays $\{(q_1,1), (q_1,-1), (q_2,1), (q_2,-1), \cdots ,(q_k,1), (q_k,-1)\}$ in $\widetilde{N}$.
 Note that $\widetilde{\beta}$ is dual to $g$ considered as a face in $B_{\Delta}$.  Moreover, we have
 \begin{eqnarray*}
&N(\widetilde{\beta})&= \frac{\nf{\widetilde{N}}{2\widetilde{N}}}{\nf{\widetilde{\beta}\cap \widetilde{N}}{\widetilde{\beta}\cap 2\widetilde{N}}} \\
& & \cong \frac{\nf{N}{2N}\oplus <e_{d+1}>}{\nf{\beta \cap N}{\beta\cap 2N} \oplus <e_{d+1}>} \\
& & \cong \frac{\nf{N}{2N}}{\nf{\beta\cap N}{\beta\cap 2N} } \\
 & & = N(\beta).
 \end{eqnarray*}
 This gives that inclusion of $[V(\widetilde{\beta})]$ in $H_{d-2}(X_{\Xi}(\mathbb{R}))$ is also represented cellularly by $C_{\beta}$.
Similarly, for each $\sigma_i \in \Sigma$ appearing in the chain $C$, $\widetilde{\sigma_i}:=\sigma_i \times [-1,1]\in \Xi$ satisfies $N(\widetilde{\sigma_i}) \cong N(\sigma_i)$.
Hence, $C$ can be viewed as a chain in $C_{d-1}(X_{\Xi}(\mathbb{R}))$.  We have
 \begin{equation}\label{tildemap}
 \widetilde{\partial} C = C_{\beta}\end{equation} where $\widetilde{\partial}$ is the cellular boundary map for
$X_{\Xi}(\mathbb{R})$.  Equation (\ref{tildemap}) holds because $\widetilde{\sigma_i} < \gamma$ in $\Xi$ implies $\gamma$ must be
of the form $\gamma=\sigma \times [-1,1]$ with $\sigma \in \Sigma$.
As $[V(\widetilde{\beta})]$ generates $A^T_{d-2}(X_{\Xi})$, the map $A^T_{d-2}(X_{\Xi})\longrightarrow H_{d-2}(X_{\Xi}(\mathbb{R}))$
must be zero.  However, the ranks of the entries $\overline{E}^1$ for $X_{\Xi}$ have the form (\ref{e1smooth}) and $\overline{E}^1_{d-2,d-2}$
is in the $4$th column from the right.  There cannot be higher boundaries with target $\overline{E}^1_{d-2,d-2}$ which contradicts
the fact that $A^T_{d-2}(X_{\Xi})\longrightarrow H_{d-2}(X_{\Xi}(\mathbb{R}))$ is the zero map.  Hence, Lemma \ref{collapse} holds and
the spectral sequence $\overline{E}^r$ for $X_{\Sigma}$ collapses at $\overline{E}^1$.
\begin{corollary}
If $\Sigma^*$ consists of $\mathbb{Z}_2$ regular cones then $X_{\Sigma}$ is maximal.
\end{corollary}

\begin{remark}
We have proved the maximality of toric varieties associated to the Fano polyhedra.
\begin{definition}
 Let $\mathrm{vert}\Delta= \{v_1, v_2, \cdots , v_n\}$.  The $d$ dimensional polytope $\Delta$ is a \emph{Fano polyhedron} provided
\begin{enumerate} \item $0 \in \mathrm{int}\Delta$
\item Each face of $\Delta$ is a simplex
\item If $v_{i_1}, v_{i_2}, \cdots ,v_{i_d}$ are the vertices of a $(d-1)$ dimensional face of $\Delta$ then $$\mathrm{det}[v_{i_1}v_{i_2} \cdots v_{i_d}]=\pm 1.$$
\end{enumerate} \end{definition}
If $\Delta$ is a Fano polyhedron then $\Delta^*$ defines one of the so called smooth toric Fano manifolds. A classification of the Fano polyhedra is known for dimension at most 4.  There are $5$ Fano polyhedra of dimension $2$.  Batyrev classified the $18$ Fano polyhedra of dimension 3 in [Bat3] and the $123$ Fano polyhedra of dimension $4$ in [Bat2].
\end{remark}
\section{Examples}\label{exsec}
In this section, we illustrate our results with two examples.  We use $\mathtt{torhom}$ [Fra] to compute the $\mathbb{Z}_2$ Hodge spaces $H_{pq}(\Sigma)$ and the $f$-vector of the polytope $\Delta$.

\begin{example}
A seven dimensional example.
\end{example}
We define $\Delta$ to be the convex hull of the following nine vertices. $$\{-e_1, -e_2, -e_3, -e_4, -e_5, - e_6, - e_7, e_1+e_2+e_3+e_4, e_5+e_6+e_7\}$$
The $f$-vector for $\Delta$ is $(9, 36,84,125 ,120,70,20)$.  The polar polytope $\Delta^*$ is the product  $P_3 \times P_4$, where $P_i$ is the $i$ dimensional simplex.  Moreover, $\Delta^*$ defines the nonsingular toric variety  $\mathbb{P}^3\times \mathbb{P}^4$.  Below are the ranks of the $\mathbb{Z}_2$ Hodge spaces $H_{pq}(\Sigma)$.
$$q
 \begin{matrix} & & & & &  & &  1\\ && &&& &15&2\\  && &&&1&31&1\\ &&&&1&0&34&0\\  &&&1&0&0&21&0 \\&&1&0&0&0&7&0 \\ &1&0&0&0&0&1&0 \\1 & 0& 0 & 0 & 0 & 0 & 0 & 0  \\ && & p && \end{matrix} $$  
We compute the $\mathbb{Z}_2$ Betti numbers for the real points $X_{\Sigma}(\mathbb{R})$ by adding along the columns $$[1\; \; \:1\; \; \:1\; \; \: 1\; \; \: 1\; \; \: 1\; \; \: 109\; \; \: 4]$$ 
and for the complex points $X_{\Sigma}(\mathbb{C})$ by adding along the diagonals.  $$[1\; \; \: 0\; \; \: 1 \; \; \: 0\; \; \: 1\; \; \: 0\; \; \: 1\; \; \: 1\; \; \: 8\; \; \: 21\; \; \: 35\; \; \:31\; \; \: 16\; \; \: 2\; \; \: 1]$$
\begin{example}
A six dimensional example.
\end{example}
 Let $\Delta$ be the convex hull of the twelve vertices below. $$\{-e_1,\pm e_2, \pm e_3, \pm e_4, \pm e_5, e_6, e_1-e_6,-e_1-e_2-e_3-e_4-e_5-e_6\}$$
The $f$-vector for $\Delta$ is $(12,62,174,267,207,64)$.  The polar polytope $\Delta^*$ has $64$ vertices.
Using $\mathtt{polymake}$ we determine that $\Delta^*$ defines a toric variety with isolated singularities.
That is, $\Sigma^*_{\leq 5}$ consists of $\mathbb{Z}_2$ regular cones.  The results of Section \ref{vansec} and Section \ref{chowsec}
give the $\mathbb{Z}_2$ Hodge spaces $H_{pq}(\Sigma)$ for $p<4$, as shown by the $\mathtt{torhom}$ computation below.
$$q
 \begin{matrix} & & & & &  &  1\\ &&&& &58&0\\  &&&&15&113&0\\ &&&1&38&96&0\\  &&1&0&34&45&0 \\&1&0&0&10&10&0 \\ 1 &  0 & 0 & 0 & 0 & 0 & 0  \\ && & p && \end{matrix} $$
We do not have the theory to guarantee collapsing of the spectral sequence $\overline{E}^r$ at $\overline{E}^1$.  However, using $\mathtt{torhom}$ we compute the $\mathbb{Z}_2$ Betti numbers for the real points $X_{\Sigma}(\mathbb{R})$.  $$[1\; \; \: 1\; \; \: 1\; \; \: 1\; \; \: 97\; \; \: 322\; \; \: 1]$$
We conclude that the spectral sequence $\overline{E}^r$ collapses at $\overline{E}^1$.  Again, we use Equation (\ref{smith}) to obtain the collapsing of
the spectral sequence $E^r$ at $E^2$.  We can therefore compute the $\mathbb{Z}_2$ Betti numbers for $X_{\Sigma}(\mathbb{C})$ by adding along the diagonals. $$[1\; \; \: 0\; \; \: 1\; \; \:0\; \; \:1\; \; \:10\; \; \: 45\; \; \: 83\; \; \: 111\; \; \: 113\; \; \: 58\; \; \: 0\; \; \: 1]$$

 \end{document}